\documentclass[%
 aip,
 amsmath,amssymb,
 reprint,%
]{revtex4-2}
\usepackage{amssymb}
\usepackage{amsthm}
\usepackage{graphicx} 
\usepackage{float}
\usepackage{color}
\usepackage[scanall]{psfrag} 
\usepackage[utf8]{inputenc}
\usepackage{amsmath}
\usepackage{amssymb}
\usepackage{mathtools}
\usepackage{bbm}
\usepackage{caption}
\usepackage{subcaption}
\usepackage{comment}
\usepackage{cool}
\usepackage{pgffor}
\usepackage{dcolumn}
\usepackage{bm}
\usepackage[T1]{fontenc}
\usepackage{mathptmx}
\usepackage{etoolbox}
\usepackage{soul}
\usepackage{enumerate}
\usepackage{hyperref}
\hypersetup{
    colorlinks=true,
    linkcolor=blue,
    filecolor=magenta,      
    urlcolor=cyan,
    pdftitle={Overleaf Example},
    pdfpagemode=FullScreen,
    }
\makeatletter
\newcommand{\@giventhatstar}[2]{#1\;\middle|\;#2)}
\newcommand{\@giventhatnostar}[3][]{#1#2\;#1|\;#3#1}
\newcommand{\giventhat}{\@ifstar\@giventhatstar\@giventhatnostar}
\makeatother
\newcommand{\FoxH}[9]{%
  H^{#1, #2}_{#3, #4} \left[ #5 \;\middle|\;
  \begin{matrix}
  \foreach \x [count=\xi] in {#6}{%
\foreach \y [count=\yi] in {#7}{%
\ifnum \xi=\yi%
\quad (\x, \y)
\fi
}
} 
  \\[0.5em]
  \foreach \x [count=\xi] in {#8}{%
\foreach \y [count=\yi] in {#9}{%
\ifnum \xi=\yi%
\quad (\x, \y)
\fi
}
}
  \end{matrix} \right]%
}
\newtheorem{theorem}{Theorem}[section]
\newtheorem{proposition}{Proposition}[section]
\theoremstyle{definition}
\newtheorem{definition}{Definition}[section]

\newtheorem{remark}{Remark}[section]

\begin{document}
\title{Riemann-Liouville fractional Brownian motion with random Hurst exponent}
\author{Hubert Woszczek}
\email{Corresponding author: hubert.woszczek@pwr.edu.pl}

    \affiliation{Faculty of Pure and Applied Mathematics, Hugo Steinhaus Center, Wroc{\l}aw University of Science and Technology, 50-370 Wrocław, Poland}%

\author{Agnieszka Wyłomańska}

    \affiliation{Faculty of Pure and Applied Mathematics, Hugo Steinhaus Center, Wroc{\l}aw University of Science and Technology, 50-370 Wrocław, Poland}%

\author{Aleksei Chechkin}
\affiliation{Faculty of Pure and Applied Mathematics, Hugo Steinhaus Center, Wroc{\l}aw University of Science and Technology, 50-370 Wrocław, Poland}%
\affiliation{Institute for Physics \& Astronomy, University of Potsdam, 14476, Potsdam-Golm, Germany}

\affiliation{German-Ukrainian Core of Excellence Max Planck Institute of Microstructure Physics, Weinberg 2, 06120 Halle (Saale), Germany}
\begin{abstract}
We examine two stochastic processes with random parameters, which in their basic versions (i.e., when the parameters are fixed) are Gaussian and display long range dependence and anomalous diffusion behavior, characterized by the Hurst exponent. Our motivation comes from biological experiments, which show that the basic models are inadequate for accurate description of the data, leading to modifications of these models in the literature through introduction of the random parameters. The first process, fractional Brownian motion with random Hurst exponent (referred to as FBMRE below) has been recently studied, while the second one, Riemann-Liouville fractional Brownian motion with random exponent (RL FBMRE) has not been explored. To advance the theory of such doubly stochastic anomalous diffusion models, we investigate the probabilistic properties of RL FBMRE and compare them to those of FBMRE. Our main focus is on the autocovariance function and the time-averaged mean squared displacement (TAMSD) of the processes. Furthermore, we analyze the second moment of the increment processes for both models, as well as their ergodicity properties. As a specific case, we consider the mixture of two point distributions of the Hurst exponent, emphasizing key differences in the characteristics of RL FBMRE and FBMRE, particularly in their asymptotic behavior. The theoretical findings presented here lay the groundwork for developing new methods to distinguish these processes and estimate their parameters from experimental data.

\end{abstract}

\maketitle
\onecolumngrid 
\begin{quotation}
One of the most well-known anomalous diffusion processes is fractional Brownian motion (FBM) which is a self-similar Gaussian process with stationary, power-law correlated increments. Another related process, Riemann-Liouville fractional Brownian motion (RL FBM), is also self-similar, Gaussian, but its increments are no longer stationary. Like FBM, the RL FBM is considered as one of the fundamental processes exhibiting anomalous diffusion properties and long range dependence. Both FBM and RL FBM share certain probabilistic characteristics with their correlation and diffusion properties determined by the self-similarity index known as the Hurst exponent $H$. However, single particle tracking experiments have uncovered highly complex phenomena in biological cells that cannot be explained within the framework of self-similar processes. Motivated by such experiments, this paper explores processes that retain the key properties of classical models like FBM and RL FBM at the level of individual trajectories, but with a Hurst exponent randomly varying from trajectory to trajectory. We examine the probabilistic properties of these modified versions of FBM and RL FBM, referred to as FBMRE and RL FBMRE, and highlight the key differences between these processes. The theoretical results presented here serve as a foundation for advanced statistical methods that can differentiate between these processes.
\end{quotation}
\section{Introduction}
In this paper we discuss two stochastic processes with random parameters. Both models, in their basic versions (with constant parameters) are Gaussian distributed and are well-known in the literature for exhibiting anomalous diffusion behavior, characterized by the non-linear time dependence of the second moment $ \mathbb{E}\left[ X^2(t)\right] \sim t^\mu$ \cite{Bouchaud1990}. 

The first process in its basic form is fractional Brownian motion (FBM). This process was originally introduced by Kolmogorov \cite{kol140} and subsequently developed by Mandelbrot and van Ness \cite{Mandelbrot1968}, particularly in the context of analyzing economic time series. FBM is characterized by the Hurst exponent $H\in(0,1)$, which governs its anomalous diffusion behavior, with $\mu=2H$. It is the only Gaussian self-similar process with stationary, power-law correlated increments and is closely linked to long-range dependence phenomena \cite{beran}. FBM has been applied in a wide range of fields, including hydrology \cite{https://doi.org/10.1029/97WR01982,BENSON2013479}, telecommunications and signal processing \cite{8246429,7448970,984735}, image analysis \cite{605414,6879494}, economics \cite{ROSTEK201330,10.1063/5.0054119,XIAO2010935}, and biological systems, such as single-particle tracking experiments \cite{weiss2012,franosch2013,metzler2014,szarek2022statistical,gleb2019}, among many others.

The second considered process, namely Riemann-Liouville fractional Brownian motion (RL FBM) was introduced by L\'evy \cite{Levy1953} in its basic form and also mentioned by Mandelbort and van Ness \cite{Mandelbrot1968}. RL FBM, similarly to FBM, is characterized by the Hurst exponent $H\in(0,1)$ responsible for the anomalous diffusion behavior. It is also self-similar, Gaussian process. However, its increments are no longer stationary and it is only self-similar, when marginal distributions are concerned \cite{Marinucci1999}. RL FBM was applied in modelling single-file diffusion \cite{Lim2009} and financial time series \cite{Picozzi2002}. 

Although both discussed processes, namely FBM and RL FBM, are widely discussed in the literature and have found many interesting applications, they seems to be inadequate for the complex data corresponding to advanced experiments like single particle tracking. Indeed, modern experimental data indicate that biological cells exhibit anomalous diffusion at the level of individual trajectories. However, the anomalous diffusion exponent varies from one trajectory to another, as demonstrated for instance in \cite{wang2018,benelli2021sub,speckner2021single,cherstvy2019non}. In such cases, a natural extension of standard anomalous diffusion models is to incorporate randomness in the parameters characterizing the process. This concept of doubly stochastic behavior forms the foundation of the superstatistics approach \cite{s_47,Beck_2005} and has been further developed in the "diffusing diffusivity" framework introduced in \cite{Chubynsky2014}, which has since been explored in various studies, as summarized in the review \cite{West2023}.

The FBM with random Hurst exponent (FBMRE) was discussed recently \cite{Balcerek2022}, revealing interesting phenomena such as accelerated diffusion and persistence transitions over time. The  extension of FBM, which allows the Hurst exponent to be a stationary random process, has been discussed in the mathematical literature \cite{levy1995multifractional,ayachetaqqu2005}.   We also direct the readers to the recent works exploring FBM-based models with different scenarios for random parameters \cite{han2020deciphering,korabel2021local,balcerek2023modelling}. There also exist the approaches allowing to distinguish between  FBM and FBMRE processes, see e.g. \cite{10.1063/5.0201436}. Besides FBM, other anomalous diffusion models have been generalized in a similar way. For example, in \cite{Woszczek2024,Santos2022}, the authors examined scaled Brownian motion with random exponent (SBMRE). Other anomalous diffusion models  with random parameters were also discussed in recent years. We mention here doubly stochastic version of continuous time random walk \cite{Arutkin2024} and L\'evy-walk-like Langevin dynamics with random parameters \cite{Chen2024}. We note that FBM with random diffusivity and constant Hurst index also have been studied \cite{Mackala2019, Wang2020, Wang2020a}.

The main goal of our research is to examine the fundamental properties of RL FBM with a random Hurst exponent (RL FBMRE) and compare them with the corresponding characteristics of FBMRE.  The focus is primarily on the autocovariance function of the discussed processes and the time-averaged mean squared displacement (TAMSD), a key statistic used to identify anomalous diffusion behavior in real data.   
We also investigate the second moment of the increment processes for both RL FBMRE and FBMRE, along with their ergodicity properties. Additionally, we explore the H{\"o}lder continuity of the processes' trajectories. These characteristics are examined for a general distribution of the Hurst exponent over the interval $(0,1)$. As a specific case, we analyze the mixture of two point distributions of the random variable $H$, highlighting the key differences in the characteristics between RL FBMRE and FBMRE, especially in their asymptotic behavior. The theoretical results presented serve as a foundation for developing new statistical techniques to distinguish these processes and estimate their parameters from experimental data.

The remainder of the paper is organized as follows. In Section \ref{preli}, we recall the basic theorems that are useful for the subsequent parts of the manuscript. In Section \ref{FBMII_sec}, we present the main characteristics of RL FBM and compare them with those of FBM. Section \ref{sec:gen}, the core of this article, discusses the main characteristics of RL FBMRE and FBMRE for a general distribution of the Hurst exponent. In Section \ref{sec:two_point}, we focus on a special case of the Hurst exponent distribution, specifically the mixture of two point distributions. Here, we also explore the asymptotic behavior of the characteristics for both processes with random Hurst exponents, highlighting the key differences between RL FBMRE and FBMRE. The final section concludes the paper. The Appendix contains the proofs.

\section{Preliminaries}\label{preli}
\noindent In this section, we introduce all  definitions and notations used in the following sections. We assume that all processes take values in $\mathbb{R}$. If we do not express it explicitly, we assume, that corresponding probability space is $\left( \Omega, \mathcal{F}, \mathbb{P}, \left\{\mathcal{F}_{t}\right\}\right)$.

\noindent We start with the general definition of the ergodic dynamical system \cite{Lasota1994}. 
\begin{definition}\label{ergodicitygeneral}
    Let $\left(Y, \mathcal{F}, \mu, S\right)$ be a measure-preserving dynamical system, where $Y$ is the state space, $\mathcal{F}$ is the $\sigma$-algebra defined on $Y$, $\mu$ is a probability measure, and $S: Y \to Y$ is a measure-preserving transformation. We call the set $A \in \mathcal{F}$ invariant if $S^{-1}\left(A\right) = A$. We call $\left(Y, \mathcal{F}, \mu, S\right)$ \textit{ergodic} if every invariant set $A \in \mathcal{F}$ is trivial, i.e. $\mu\left(A\right) = 0$ or $\mu\left(X \backslash A\right) = 0$.  
\end{definition}
\noindent For our purpose we need to translate the general Definition \ref{ergodicitygeneral} to the language of stochastic processes. Let $\left\{ X\left(t\right)\right\}_{t \in \mathbb{N}_0}$ be the stationary stochastic process. The assumption of stationarity of the process is necessary because it guarantees the invariance of the measure preserving transformation, which is usually time shift in the case of stochastic processes. The stochastic process $\left\{ X\left(t\right)\right\}_{t \in \mathbb{N}_0}$ in its canonical representation can then be identified with its law, which is the probability measure on the space $\mathbb{R}^{\mathbb{N}_0}$. Thus, we identify space $\mathbb{R}^{\mathbb{N}_0}$ with state space and probability measure $\mathbb{P}$ with probability measure $\mu$. The $\sigma$-algebra $\mathcal{F}$ is then the canonical $\sigma$-algebra $\mathfrak{B}$ generated by the cylinder sets. Transformation $S$ in this case is (time) shift transformation, namely $S_t f\left(s\right) = f\left(s+t\right)$. Finally, the dynamical system has the following form $\left(\mathbb{R}^{\mathbb{N}_0}, \mathfrak{B}, \mathbb{P}, S_t\right)$.

\noindent Next, we recall the convinient theorem, which gives necessary and sufficient condition for the stochastic process to be ergodic in terms of characteristic function. This theorem will be crucial in the proof of the lack of ergodicity of FBMRE.
\begin{theorem}\label{parzen}
    Let $\left\{X\left(t\right)\right\}_{t \in \mathbb{N}}$ be strictly stationary stochastic process. Then $\left\{X\left(t\right)\right\}_{t \in \mathbb{N}}$ is ergodic iff. the following condition holds \cite{Parzen1958}
    \begin{equation}
        \begin{aligned}
            \lim_{n \to  \infty} \frac{1}{n+1} \sum_{\tau=0}^n \mathbb{E}\left[\exp\left(i \sum_{j=1}^K z_j \left(X\left(t_j + h \right) - X\left(t_j + h + \tau\right)\right)\right)\right]  = \\ = \left|\mathbb{E}\left[\exp\left(i \sum_{j=1}^K z_k X\left(t_j+h\right)\right)\right]\right|^2.
        \end{aligned}
        \end{equation}
\end{theorem}
\noindent Next, we recall Jensen's inequlaity, which is also important for the proof of nonergodicity of FBMRE.
\begin{theorem}
    Let $\mu$  be positive-defined measure on some $\sigma$-algebra of the set $\Omega$ and $\mu\left(\Omega\right)=1$. If a real function $f \in L^1\left(\mu\right)$, $a < f\left(x\right)<b$ for all $x\in \Omega$ and if $\varphi$ is a convex function on the interval $\left(a,b\right)$, then \cite{Rudin}
    \begin{equation}\label{jensenineq}
        \varphi\left(\int_{\Omega}f d\mu\right) \leq \int_{\Omega}\left(\varphi \circ f \right)d\mu.
    \end{equation}
If $\varphi$ is concave the reverse inequality is true.
\end{theorem}
\noindent In the context of probability theory for some random variable $X$ inequality \eqref{jensenineq} takes the following form
\begin{equation}\label{jensenprob}
    \varphi\left(\mathbb{E}\left[X\right]\right) \leq \mathbb{E}\left[\varphi\left(X\right)\right].
\end{equation}
\noindent Finally we recall the Kolmogorov continuity theorem, which is convenient tool to show the regularity of trajectories of FBMRE and RL FBMRE \cite{Schilling2012}.
\begin{theorem}\label{kolmogorov}
    Let $\left\{X\left(t\right)\right\}_{t \in \left[0,1\right]}$ be the real-valued stochastic process. Suppose that there exists positve constants q, c, $\epsilon$ such that for every $s,t \in \left[0, 1\right]$ 
    \begin{equation}
        \mathbb{E}\left[\left|X\left(t\right) - X\left(s\right)\right|^q\right] \leq c\left|t-s\right|^{1+\epsilon}.
    \end{equation}
    Then there exists a modification $\left\{\tilde{X}\left(t\right)\right\}_{t \in \left[0,1\right]}$ of $\left\{X\left(t\right)\right\}_{t \in \left[0,1\right]}$ which is almost surely H{\"o}lder continuous with exponent $\kappa \in \left(0, \epsilon/q\right)$.
\end{theorem}
\section{Riemann-Liouville fractional Brownian motion}\label{FBMII_sec}
In this section, we recall the definition of the RL FBM model and present its main characteristics. Additionally, we highlight the key differences between RL FBM and FBM.
\begin{definition}\label{fbmiidef}
The Riemann-Liouville fractional Brownian motion  $\left\{B_H^*\left(t\right)\right\}_{t \geq 0}$ is defined as follows \cite{Levy1953}
\begin{equation}
    B_H^*\left(t\right) = \int_0^t \sqrt{2H} \left(t-s\right)^{H - \frac{1}{2}} dB\left(s\right),
\end{equation}
where $H \in \left(0, 1\right)$ is the Hurst exponent and $\left\{B\left(t\right)\right\}_{t \geq 0}$ is the ordinary Brownian motion. The prefactor $\sqrt{2H}$ guarantees the fullfillment of
\end{definition}
\begin{proposition}\label{varfbmii}
The second moment of RL FBM is given by \cite{Lim2001}
\begin{equation}\label{secon_FBM}
\mathbb{E}\left[\left(B_H^*\left(t\right)\right)^2\right] = t^{2H}.
\end{equation}
\end{proposition}

\noindent In mathematical literature the RL FBM is also called FBM type II process \cite{Marinucci1999}. 

\noindent The RL FBM is a zero-mean, Gaussian, self-similar process with the probability density function (pdf) given by
\begin{equation}\label{pdf_FBM}
    p\left(x, t \right) = \frac{1}{\sqrt{2\pi t^{2H}}} \exp\left\{-\frac{x^2}{2t^{2H}}\right\}, \quad x \in \mathbb{R}.
\end{equation}

\noindent We recall, the FBM $\left\{B_H\left(t\right)\right\}_{t \geq 0}$is also the zero-mean Gaussian process with the self-similarity property. Additionally, its pdf and the second moment are exactly the same as for RL FBM, namely they are given in Eq. (\ref{pdf_FBM}) and (\ref{secon_FBM}), respectively. 

\noindent Next we recall the formula for the autocovariance function of RL FBM.
\begin{proposition}\label{acvffbmii}
   The autocovariance function of RL FBM is given by the following formula \cite{Lim2001}
    \begin{equation}\label{acvffbmiihipergeo}
        Cov\left(B_H^*\left(t\right), B_H^*\left(t+\tau\right)\right) = \frac{2H \left(t+\tau\right)^{H - \frac{1}{2}}t^{H + \frac{1}{2}}}{H + 1/2} \Hypergeometric{2}{1}{\frac{1}{2} - H, 1}{\frac{3}{2} + H}{\frac{t}{t + \tau}},
    \end{equation}
where $s<t$ and $\Hypergeometric{2}{1}{a,b}{c}{z}$ is Gauss hypergeometric function. Equivalently, the formula  (\ref{acvffbmiihipergeo}) can be written as \cite{Marinucci1999}
\begin{equation}\label{acvffbmiialternative}
Cov\left(B_H^*\left(t\right), B_H^*\left(t+\tau\right)\right) = \frac{1}{2}\left( \left(t + \tau\right)^{2H} + t^{2H} - \mathbb{E}\left[\left(b_{H}^{*\tau}\left(t\right)\right)^2\right]\right).
\end{equation}
In the above formula  $\{b_{H}^{*\tau}\left(t\right)\}_{t\geq 0}$ is the  increment of RL FBM,
\begin{eqnarray}\label{inc_FBMII}b_{H}^{*\tau}\left(t\right) = B_H^*\left(t + \tau\right) - B_H^*\left(t\right).\end{eqnarray}
\end{proposition}
\noindent We recall, for FBM model, the autocovariance function has the following form (assuming the diffusion coefficient is equal to one)\begin{equation}\label{FBM1}
Cov\left(B_{{H}}\left(t\right), B_{{H}}\left(t+\tau\right)\right) = \frac{1}{2}\left(\left|t+\tau\right|^{2{H}} + \left|t\right|^{2{H}} - \left|\tau\right|^{2{H}}\right).
    \end{equation}  
\noindent Let us note, the formula (\ref{FBM1}) can also be represented as in Eq. (\ref{acvffbmiialternative}), however in this case the second moment of the increment process of RL FBM is replaced by the second moment of the increments of FBM defined similar as in Eq. (\ref{inc_FBMII}).

\noindent In the next proposition we present the asymptotics for the autocovariance function of the RL FBM for long and short times.
\begin{proposition}\label{fbmiicovasympt}
The  autocovariance function of RL FBM has the following asymptotic
  \begin{equation}\label{fbmiicovasympteq}
Cov\left(B_H^*\left(t\right), B_H^*\left(t+\tau\right)\right) \sim
    \begin{cases}
        \tau^{2H}\left(\frac{t}{\tau}\right)^{H+\frac{1}{2}} \left(\frac{2H}{\frac{1}{2} + H} - \frac{2H\left(\frac{1}{2}-H\right)}{\frac{3}{2} + H}\frac{t}{\tau} +  \frac{H\left(\frac{1}{2}-H\right)\left(\frac{3}{2}-H\right)}{\frac{5}{2} + H}\left(\frac{t}{\tau}\right)^{2}\right), ~~ \text{ for }~~ t/\tau << 1,\\ 
         \frac{1}{2}t^{2H}\left(2 + 2H\frac{\tau}{t} - \frac{2H\Gamma\left(H + \frac{1}{2}\right)^2}{\Gamma\left(1+2H\right)\sin\left(\pi H\right)}\left(\frac{\tau}{t}\right)^{2H} + H\left(2H-1\right)\left(\frac{\tau}{t}\right)^2\right), ~~\text{ for } ~~t/\tau >> 1.
    \end{cases}
\end{equation} 
\begin{proof}
\noindent In case $t/\tau<<1$ we use the representation of autocovariance of RL FBM via hypergeometric function \eqref{acvffbmiihipergeo} and asymptotic of $I_3$ from the  proof of Proposition \ref{incvarfbmiiasympt} (see Appendix \ref{AppA}).
\\
\noindent In case when $t/\tau>>1$ we apply the Taylor expansion for $\left(t+\tau\right)^{2H}$ in equation \eqref{acvffbmiialternative} and formula \eqref{structfuncfbmiiasympt} presented below.
\end{proof}
\end{proposition}
\noindent We recall that the notation  $f\left(x\right) \sim g\left(x\right)$ for two real functions means the following 
\begin{equation}
    \lim_{x \to x_0} \frac{f\left(x\right)}{g\left(x\right)} = 1,
\end{equation}
where $x_0$ (possibly "equal" to $\infty$) depends on the context and can be easily read from it.
For $a,b \in \mathbb{R}$, when we write $a<<b$ we mean, that "b is sufficiently larger than a".\\
 
\begin{proposition}
     The autocovariance function of  FBM behaves asymptotically as
\begin{equation}\label{covfbmtpasympt}
    Cov\left(B_H\left(t\right), B_H\left(t+\tau\right)\right) \sim
    \begin{cases}
        \frac{1}{2}\tau^{2H}\left(2H\frac{t}{\tau}  + \left(\frac{t}{\tau}\right)^{2H} + H\left(2H-1\right)\left(\frac{t}{\tau}\right)^2 \right), ~~\text{ for }~~ t/\tau<<1\\
        \frac{1}{2}t^{2H}\left(2 + 2H\frac{\tau}{t}  - \left(\frac{\tau}{t}\right)^{2H} + H\left(2H-1\right)\left(\frac{\tau}{t}\right)^2\right), ~~\text{ for }~~ t/\tau>>1.
    \end{cases}
\end{equation}
\begin{proof}
In order to obtain the asymptotics of the autocovariance function for FBM we expand the term $\left(t + \tau\right)^{2H}$ in the Taylor series and consider its first three components.
\end{proof}
\end{proposition}
\noindent Let us note that for $H=1/2$ the autocovariance functions of RL FBM and FBM are the same and equal $t$, as it should be. It can be easily seen from equations \eqref{acvffbmiihipergeo}-\eqref{FBM1}. However, for $H\neq1/2$ we observe that the orders of the dominant terms for short times are different. Specifically, in the former case, the dominant term for FBM is of the order $\mathcal{O}\left(\tau^{2H}\frac{t}{\tau}\right)$, whereas for RL FBM, the dominant term is of the order $\mathcal{O}\left(\tau^{2H}\left(\frac{t}{\tau}\right)^{H+\frac{1}{2}}\right)$. For $t/\tau>>1$ in both cases the dominant term is $t^{2H}$.

\noindent Now we go to the expectation of the TAMSD of RL FBM, which is very important quantity measured in experiments. We recall the TAMSD of the process $\left\{ X\left(t\right)\right\}_{t\geq 0}$ is defined as \cite{metzler2014}
    \begin{equation}\label{tamsd}
    \delta\left(\tau\right) = \frac{1}{T - \tau} \int_0^{T-\tau} \left(X\left(t + \tau\right) - X\left(t\right)\right)^2 dt,
\end{equation}
where $\tau \in \left[0, T\right)$ is the time lag of measured time series (the width of a sliding window), and $T \in \left(0, \infty\right)$ is a time horizon (the trajectory length).
\begin{proposition}\label{th1}   The expectation of TAMSD of RL FBM has the form
    \begin{equation}\label{fbmiitamsd}
        \begin{aligned}
            \mathbb{E}\left[\delta_H^*\left(\tau\right)\right] = 
            \frac{1}{2H+1}\frac{T^{2H+1} - \tau^{2H+1}}{T- \tau} + \frac{1}{2H+1}\left(T - \tau\right)^{2H} - \\ - \frac{4H}{\left(H+\frac{1}{2}\right)\left(H+\frac{3}{2}\right)} \frac{\left(T - \tau\right)^{H + 1/2}}{\tau^{1/2 - H}}\Hypergeometric{2}{1}{\frac{1}{2} + H, \frac{1}{2}-H}{\frac{5}{2} + H}{\frac{\tau-T}{\tau}}.
    \end{aligned}
    \end{equation}
    The proof of this proposition is presented in Appendix \ref{AppA}. When we put $H=\frac{1}{2}$, we obtain $\mathbb{E}\left[\delta_{\frac{1}{2}}^*\left(\tau\right)\right] = \tau$ as it should be.
    \end{proposition}
\noindent Let us note that for FBM the expectation of the TAMSD is the same as the second moment of the corresponding increment process, namely $\mathbb{E}\left[\delta_H\left(\tau\right)\right] =\tau^{2H}$.\\
\begin{proposition}\label{tamsdfbmiiasympt}
    For $T/\tau >>1$ the asymptotic of the expectation of TAMSD of RL FBM reads
    \begin{equation}\label{tamsdfbmiiasympteq}
        \begin{aligned}
            \mathbb{E}\left[\delta_H^*\left(\tau\right)\right] \sim
         \frac{2H\Gamma\left(\frac{1}{2}+H\right)^2\tau^{2H}}{\Gamma\left(1+2H\right)\sin\left(\pi H\right)}. 
        \end{aligned}
    \end{equation}
\noindent The proof of this proposition is presented in Appendix \ref{AppA}.
    \end{proposition}
    \noindent It is worth to highlight that for $H=1/2$, asymptotics of the expectation of TAMSD for RL FBM and FBM are the same, however for $H\neq 1/2$, they differ (the prefactors are different).\\
\noindent As the last characteristic, we investigate the second moment of the increment process of RL FBM  defined in (\ref{inc_FBMII}). Let us note, in contrast to FBM, the increment process of RL FBM is not stationary.
\begin{proposition}\label{incvarfbmii}
    The second moment of the increment process of RL FBM is given by the following formula \cite{Lim2001}
    \begin{equation}\label{incvarfbmiieq}
        \mathbb{E}\left[\left(b_{H}^{*\tau}\left(t\right)\right)^2 \right] = 2H\tau^{2H}\left(I\left(t, \tau; H\right) + \frac{1}{2H}\right),
    \end{equation}
    where \begin{eqnarray}\label{II}I\left(t, \tau; H\right) = \int_0^{t/\tau} \left[\left(1 + u\right)^{H- 1/2} - u^{H - 1/2}\right]^2 du.\end{eqnarray}
\end{proposition}
\noindent To derive asymptotic of the second moment of the increment process for $t/\tau<<1$ in the next proposition we present the equivalent representation of the discussed characteristics.
\begin{proposition}\label{incrementscovfbmiialter}
    If $H \neq 1/2$, then the second moment of the increment process of RL FBM can be expressed as
    \begin{equation}\label{varincfbmiifinal}
        \begin{aligned}
            \mathbb{E}\left[\left(b_{H}^{*\tau}\left(t\right)\right)^2 \right] = \left(t+\tau\right)^{2H} + t^{2H} - \frac{4H \left(t+\tau\right)^{H - \frac{1}{2}}t^{H + \frac{1}{2}}}{H + 1/2} \Hypergeometric{2}{1}{\frac{1}{2} - H, 1}{\frac{3}{2} + H}{\frac{t}{t + \tau}}.
        \end{aligned}
    \end{equation}
\begin{proof} Using the definition of the increment process of RL FBM and applying Propositions \ref{varfbmii} and \ref{acvffbmii}, we get Eq. \eqref{varincfbmiifinal}.
    \end{proof}
    \noindent When we put $H=1/2$, then we $\mathbb{E}\left[\left(b_{H}^{*\tau}\left(t\right)\right)^2 \right] = \tau$ as it should be.
 \end{proposition}
\noindent Let us note that the second moment of the increment process of FBM is independent on $t$ and has simpler form as those presented in Eq. (\ref{incvarfbmiieq}) or Eq. (\ref{varincfbmiifinal}) for RL FBM. Namely, we have \cite{Mandelbrot1968}
\begin{eqnarray}
   \mathbb{E}\left[\left(b_{H}^{\tau}\left(t\right)\right)^2 \right] = \tau^{2H}.
\end{eqnarray}
In the next proposition we show asymptotics of the second moment of the increment proces of RL FBM.
\begin{proposition}\label{incvarfbmiiasympt}
    The second moment of the increment process of RL FBM have the following asymptotics
\begin{equation}\label{structfuncfbmiiasympt}
\mathbb{E}\left[\left(b_{H}^{*\tau}\left(t\right)\right)^2 \right] \sim
    \begin{cases}
         \tau^{2H},  \,~~\text{ for }~~ t/\tau << 1,\\
         \tau^{2H} \frac{2H\Gamma\left(H + \frac{1}{2}\right)^2}{\Gamma\left(1+2H\right)\sin\left(\pi H\right)},   \,~~\text{ for } ~~t/\tau >> 1.
    \end{cases}
\end{equation}
The proof of this proposition is presented in Appendix \ref{AppA}.
\end{proposition}

\section{Fractional Brownian motion and Riemann-Liouville fractional Brownian motion with random Hurst exponent - general case }\label{sec:gen}
In this section, we first recall the definitions of the two processes under consideration, namely FBMRE and RL FBMRE, assuming the Hurst exponent follows a general distribution over the interval $(0,1)$. We then present their main properties, paying attention to the characteristics shared by both processes, as well as those that distinguish them.
\subsection{FBMRE}
\begin{definition}\label{fbmredef11}
Let $\mathcal{H}$ be some positive random variable with values from the interval $\left(0, 1\right)$ with pdf $f_{\mathcal{H}}\left(\cdot\right)$ and a moment generating function $M_\mathcal{H}\left(\cdot\right)$,
\begin{equation}\label{moment_gen}
    M_\mathcal{H}(s) = \mathbb{E}\left[e^{s\mathcal{H}}\right] =  \int_0^1 e^{sh} f_\mathcal{H}\left(h\right)dh.
\end{equation}

\noindent Then, the fractional Brownian motion with random Hurst exponent is defined as \cite{Balcerek2022}
\begin{equation}\label{fbmredef}
    B_{\mathcal{H}}\left(t\right) = A_{\mathcal{H}}\left[\int_0^t \left(t - u\right)^{\mathcal{H} - 1/2} d\tilde{B}\left(u\right) + \int_{-\infty}^0 \left(\left(t-u\right)^{\mathcal{H} - 1/2} - \left(-u\right)^{\mathcal{H}-1/2}\right)d\tilde{B}\left(u\right) \right],
\end{equation}
where $\left\{\tilde{B}\left(t\right)\right\}_{t \in \mathbb{R}}$ is extension of the standard Brownian motion to the negative time axis,
\begin{equation}
    \tilde{B}\left(t\right) = 
    \begin{cases}
        B_1\left(t\right), \, ~~\text{ for }~~ t>0 \\ 
        B_2\left(t\right), \, ~~\text{ for }~~ t\leq0,
    \end{cases}
\end{equation}
where $\left\{B_1\left(t\right)\right\}_{t \geq 0}$ and $\left\{B_2\left(t\right)\right\}_{t \geq 0}$ are two independent standard Brownian motions, and \\ $A_{\mathcal{H}} = \sqrt{\Gamma\left(2\mathcal{H}+1\right)\sin\left(\pi \mathcal{H}\right)}/\Gamma\left(\mathcal{H}+1/2\right)$. The form of $A_{\mathcal{H}}$ ensures that under the condition $\mathcal{H}=H$ for given $ t\geq 0$, the random variable $B_{\mathcal{H}}\left(t\right)$ is zero-mean Gaussian distributed with  $\mathbb{E}\left[B_{\mathcal{H}}\left(t\right)\right] = t^{2H}$. In the above definition, we assume that $\left\{\tilde{B}\left(t\right)\right\}_{t\geq0}$ and $\mathcal{H}$ are independent. From the integral representation \eqref{fbmredef}, we immediately observe that the increments of FBMRE are strictly stationary, just as they are for FBM.
\end{definition}
\noindent In the next two propositions we recall the formulas for the pdf and q-th order absolute moments of the FBMRE.
\begin{proposition}
    The pdf of FBMRE is given by the following formula \cite{Balcerek2022}
    \begin{equation}\label{pdffbmreeq}
        p_{B_{\mathcal{H}}}\left(x, t\right) = \int_0^1 \frac{1}{\sqrt{2\pi t^{2h}}}exp\left\{\frac{-x^2}{2t^{2h}}\right\} f_\mathcal{H}\left(h\right) dh.
    \end{equation}
\end{proposition}
\begin{proposition}\label{msdfbmre}
The moment of order q of the absolute value of $\left\{B_{\mathcal{H}}\left(t\right)\right\}_{t\geq0}$ is given by \cite{Balcerek2022}
\begin{equation}
\label{q-momentfbmre}
    \mathbb{E}\left[\left|B_{\mathcal{H}}\left(t\right)\right|^q\right]  = c_q M_\mathcal{H}\left(q\log t\right),
\end{equation}
where $c_q = \frac{2^{q/2} \Gamma\left(\frac{q+1}{2}\right)}{\sqrt{\pi}}$.
\end{proposition}
\noindent In the next proposition, we present the general formula for the autocovariance function of the FBMRE, which is expressed in terms of the moment-generating function of the random variable $\mathcal{H}$ given in \eqref{moment_gen}. 
\begin{proposition}\label{acvffbmre}  
The autocovariance function of FBMRE reads\begin{equation}  
    Cov\left(B_{\mathcal{H}}\left(t\right), B_{\mathcal{H}}\left(t+\tau\right)\right) = \frac{1}{2}\left[M_\mathcal{H}\left(2\log \left|t+\tau\right|\right) + M_\mathcal{H}\left(2\log \left|t\right|\right) - M_\mathcal{H}\left(2\log \left|\tau\right|\right)\right].
\end{equation}
\begin{proof}
To obtain the autocovariance function of the considered process, we apply the law of total expectation. Thus, we have
    \begin{equation}
    \begin{aligned}
        Cov\left(B_{\mathcal{H}}\left(t\right)B_{\mathcal{H}}\left(t+\tau\right)\right) = \mathbb{E}\left[B_{\mathcal{H}}\left(t\right)B_{\mathcal{H}}\left(t+\tau\right)\right] = \mathbb{E}\left[\mathbb{E}\left[\giventhat{B_{\mathcal{H}}\left(t\right)B_{\mathcal{H}}\left(t+\tau\right)}{\mathcal{H}}\right]\right] = \\ =\mathbb{E}\left[\frac{1}{2}\left(\left|t+\tau\right|^{2\mathcal{H}} + \left|t\right|^{2\mathcal{H}} - \left|\tau\right|^{2\mathcal{H}}\right)\right]
        = \frac{1}{2}\left[M_\mathcal{H}\left(2\log \left|t+\tau\right|\right) + M_\mathcal{H}\left(2\log \left|t\right|\right) - M_\mathcal{H}\left(2\log \left|\tau\right|\right)\right].
    \end{aligned}
    \end{equation}
\end{proof}
\end{proposition}
\begin{proposition}
    The expectation of TAMSD for FBMRE is given by
   \begin{equation}\label{tamsdfbmre}
       \mathbb{E}\left[\delta_\mathcal{H}\left(\tau\right)\right] = M_\mathcal{H}\left(2\log \tau\right).
   \end{equation}
    \begin{proof}
To obtain the expectation of TAMSD, we apply Propositions \ref{msdfbmre}, \ref{acvffbmre} in the formula  \eqref{tamsd}. Thus, we obtain 
    \begin{equation}
    \begin{aligned}
    \mathbb{E}\left[\delta_\mathcal{H}\left(\tau\right)\right] = \mathbb{E}\left[\frac{1}{T - \tau} \int_0^{T-\tau} \left(B_{\mathcal{H}}\left(t + \tau\right) - B_{\mathcal{H}}\left(t\right)\right)^2 dt\right] 
    = \\ = \mathbb{E}\left[\frac{1}{T - \tau} \int_0^{T-\tau} (B^2_{\mathcal{H}}\left(t + \tau\right) - 2B_{\mathcal{H}}\left(t + \tau\right)B_{\mathcal{H}}\left(t\right) + B^2_{\mathcal{H}}\left(t\right) dt\right]
    \end{aligned}
    \end{equation}
Applying Proposition \ref{acvffbmre} and the Fubinis theorem \cite{Rudin} we obtain Eq. \eqref{tamsdfbmre}.
    \end{proof}
\end{proposition}
\noindent Next, we investigate properties of increments of FBMRE. The increment process  $\{b_{\mathcal{H}}^{\tau}\left(t\right)\}_{t\geq 0}$ of FBMRE is defined as
\begin{equation}
    b_{\mathcal{H}}^{\tau}\left(t\right) = B_{\mathcal{H}}\left(t + \tau\right) - B_{\mathcal{H}}\left(t\right).
\end{equation}
Let us recall, the increment process of FBMRE is stationary. In the next Proposition we prove the formula for its second moment.
\begin{proposition}
The second moment of increment process of FBMRE is given by
\begin{equation}
    \mathbb{E}\left[\left(b_{\mathcal{H}}^{\tau}\left(t\right)\right)^2\right] = M_\mathcal{H}\left(2\log \left|\tau\right|\right).
\end{equation}
\begin{proof}
To obtain the second moment of the increment process  we apply Propositions \ref{msdfbmre}, \ref{acvffbmre}. Thus, we have
\begin{equation}
    \begin{aligned}
        \mathbb{E}\left[\left(b_{\mathcal{H}}^{\tau}\left(t\right)\right)^2\right] = \mathbb{E}\left[\left(B_{\mathcal{H}}\left(t + \tau\right)\right)^2\right] - 2\mathbb{E}\left[B_{\mathcal{H}}\left(t + \tau\right)B_{\mathcal{H}}\left(t\right)\right] + \mathbb{E}\left[\left(B_{\mathcal{H}}\left(t\right)\right)^2\right] = M_\mathcal{H}\left(2\log \left|t+\tau\right|\right) - \\ -M_\mathcal{H}\left(2\log \left|t+\tau\right|\right) - M_\mathcal{H}\left(2\log \left|t\right|\right) +M_\mathcal{H}\left(2\log \left|\tau\right|\right) + M_\mathcal{H}\left(2\log \left|t\right|\right) = M_\mathcal{H}\left(2\log \left|\tau\right|\right).
    \end{aligned}
\end{equation}
\end{proof}
\end{proposition}
\noindent In the next theorem, which is one of the main result of this paper, we prove that the increment process of FBMRE while being stationary process, is not ergodic, even though the second moment of the increment process is equal to the expectation of the TAMSD of FBMRE.
\begin{theorem}
    The sequence of increments of FBMRE is not ergodic. 
    \begin{proof}
        We start the proof with recalling a well-known fact that sequence of increments of FBM is ergodic. The ergodicity of sequence of increments of FBM is straightforward consequence of classical  It{\^o}s result, which states that the stationary Gaussian process is mixing iff its autocovariance decays to zero at infinity. Mixing is stronger condition than ergodicity and implies it (Khinchin theorem) \cite{khinchin1949}; more information about mixing can be found, i.e., in \cite{Lasota1994}. Thus, for increments of FBM convergence of the multidimensional characteristic functions holds. By applying Theorem \ref{parzen} and Jensen's inequality \eqref{jensenprob} to the increment process of FBMRE, we obtain
        \begin{equation}
        \begin{aligned}
            \lim_{n \to  \infty} \frac{1}{n+1} \sum_{\tau=0}^n \mathbb{E}\left[\exp\left(i \sum_{j=1}^K z_j \left(b_{\mathcal{H}}^1\left(t_j + h \right) - b_{\mathcal{H}}^1\left(t_j + h + \tau\right)\right)\right)\right] = \\ 
            = \lim_{n \to  \infty} \frac{1}{n+1} \sum_{\tau=0}^n \mathbb{E}\left[\mathbb{E}\left[\giventhat{\exp\left(i \sum_{j=1}^K z_j \left(b_{\mathcal{H}}^1\left(t_j + h \right) - b_{\mathcal{H}}^1\left(t_j + h + \tau\right)\right)\right)}{\mathcal{H}}\right]\right] = \\
            = \mathbb{E}\left[\lim_{n \to  \infty} \frac{1}{n+1} \sum_{\tau=0}^n \mathbb{E}\left[\giventhat{\exp\left(i \sum_{j=1}^K z_j \left(b_{\mathcal{H}}^1\left(t_j + h \right) - b_{\mathcal{H}}^1\left(t_j + h + \tau\right)\right)\right)}{\mathcal{H}}\right]\right] = \\
            = \mathbb{E}\left[\left|\mathbb{E}\left[\giventhat{\exp\left(i \sum_{j=1}^K z_j b_{\mathcal{H}}^1\left(j+h\right)\right)}{\mathcal{H}}\right]\right|^2\right] \geq \left|\mathbb{E}\left[\exp\left(i \sum_{j=1}^K z_j b_{\mathcal{H}}^1\left(t_j+h\right)\right)\right]\right|^2.
        \end{aligned}
        \end{equation}
      The change in the order of expectation and limit is justified by the dominated convergence theorem (we refer to the monograph \cite{Rudin}), as the conditional expectation is a measurable function and the absolute value of the characteristic function is bounded by 1. The last inequality is a direct consequence of Jensen's inequality \eqref{jensenprob}. Equality holds if and only if $\mathcal{H}$ is constant, in which case FBMRE reduces to FBM. 
    \end{proof}
\end{theorem}
\noindent In the next theorem we prove that trajectories of FBMRE are almost sure H{\"o}lder continuous. 
\begin{theorem}\label{fbmrecontinuity}
Let $\left\{B_{\mathcal{H}}\left(t\right)\right\}_{t\in\left[0,1\right]}$ be the FBMRE with $\mathcal{H}$ taking values from the interval $\left(H_1, 1\right)$, $0<H_1<1$. Then, there exists a version of FBMRE such that its trajectories are almost sure H{\"o}lder continuous for all exponents $\kappa <H_1$.
\begin{proof}
In the proof we apply the Kolmogorov criterion (Theorem \ref{kolmogorov}). For FBMRE we have
    \begin{equation}
        \mathbb{E}\left[\left|B_{\mathcal{H}}\left(t\right) - B_{\mathcal{H}}\left(s\right)\right|^q\right] = c_q\int_{H_1}^{1}\left|t-s\right|^{qh}f_{\mathcal{H}}\left(h\right) dh \leq c_q\left|t-s\right|^{qH_1},
    \end{equation}
    where $c_q = \frac{2^{q/2} \Gamma\left(\frac{q+1}{2}\right)}{\sqrt{\pi}}$. We set $\epsilon = qH_1 -1 = R$ and $q = \frac{\left(1+R\right)}{H_1}$ for fixed $R$. Then we have
    \begin{equation}
        \kappa < \frac{\epsilon}{q} = \frac{RH_1}{R+1}.
    \end{equation}
    Finally, since $R$ is arbitrary, FBMRE is H{\"o}lder continuous for all exponents $\kappa <H_1$.
\end{proof}
\end{theorem}

\subsection{RL FBMRE}
\begin{definition}
Let  $\mathcal{H}$ be the random variable defined as in Definition \ref{fbmredef11}. Then, the Riemann-Liouville fractional Brownian motion with random Hurst exponent  is defined as
\begin{equation}
    B_{\mathcal{H}}^*\left(t\right) = \int_0^t \sqrt{2\mathcal{H}} \left(t-s\right)^{\mathcal{H} - \frac{1}{2}} dB\left(s\right),
\end{equation}
We assume that standard Brownian motion $\left\{B\left(t\right)\right\}_{t\geq0}$ and $\mathcal{H}$ are independent. The prefactor $\sqrt{2\mathcal{H}}$ guarantees that $\mathbb{E}\left[(B_{\mathcal{H}}^*\left(t\right))^2\right]=t^{2H}$ for each $t\geq 0$ under the condition $\mathcal{H}=H$. 
\end{definition}
\noindent Let us note that, after simple calculations, one may show that the pdf and the q-th absolute moment of RL FBMRE are the same as those for FBMRE; see formulas \eqref{pdffbmreeq} and \eqref{q-momentfbmre}, respectively. In the next proposition, we derive the autocovariance function of RL FBMRE. We note that this characteristic for RL FBMRE differs from that obtained for FBMRE; see Proposition \ref{acvffbmre}. 
\begin{proposition}\label{acvffbmreii}
The autocovariance function of RL FBMRE is given by
\begin{equation}
    Cov\left(B_{\mathcal{H}}^*\left(t\right), B_{\mathcal{H}}^*\left(t+\tau\right)\right) = \sqrt{\frac{t}{t + \tau}}\int_0^1 \frac{2h\left(t\left(t+\tau\right)\right)^h}{h+ 1/2} \Hypergeometric{2}{1}{\frac{1}{2} - h, 1}{\frac{3}{2} + h}{\frac{t}{t+\tau}}f_\mathcal{H}\left(h\right) dh,
\end{equation}
where $\prescript{}{2}{F_1 \left(\cdot, \cdot, \cdot, \cdot\right)}$ is Gauss hypergeometric function (see \cite{Beals2010} for details).
    \begin{proof}
Applying the law of total expectation and using the Proposition \ref{acvffbmii} for autocovariance of RL FBM, we have
    \begin{align*}
        Cov\left(B_{\mathcal{H}}^*\left(t\right), B_{\mathcal{H}}^*\left(t+\tau\right)\right) = \mathbb{E}\left[B_{\mathcal{H}}^*\left(t\right) B_{\mathcal{H}}^*\left(t+\tau\right)\right] = \mathbb{E}\left[\mathbb{E}\left[\giventhat{B_{\mathcal{H}}^*\left(t\right)B_{\mathcal{H}}^*\left(t+\tau\right)}{\mathcal{H}}\right]\right] = \\
        = \mathbb{E}\left[\frac{2\mathcal{H}\left(t+\tau\right)^{\mathcal{H} - \frac{1}{2}}t^{\mathcal{H} + \frac{1}{2}}}{\mathcal{H}+1/2} \Hypergeometric{2}{1}{\frac{1}{2} - \mathcal{H}, 1}{\frac{3}{2} + \mathcal{H}}{\frac{t}{t+\tau}}\right] = \\ = \sqrt{\frac{t}{t+\tau}}\int_0^1 \frac{2h\left(t\left(t+\tau\right)\right)^h}{h+ 1/2} \Hypergeometric{2}{1}{\frac{1}{2} - h, 1}{\frac{3}{2} + h}{\frac{t}{t+\tau}}f_\mathcal{H}\left(h\right) dh.
    \end{align*}
    \end{proof}
  
\end{proposition}
\noindent In the next proposition we derive the expectation of TAMSD of RL FBMRE. Like the autocovariance function, this characteristics has different form for the considered processes in the general case of $\mathcal{H}$ distribution. 
\begin{proposition}
   The expectation of TAMSD for RL FBMRE  is given by
   \begin{equation}\label{tamsdsbmre}
   \begin{aligned}
       \mathbb{E}\left[\delta_\mathcal{H}^*\left(\tau\right)\right] = \frac{1}{T - \tau} \int_0^{T-\tau} \left[M_\mathcal{H}\left(2\log \left( t + \tau\right)\right)  +  M_\mathcal{H}\left(2\log t\right)\right]dt - \\ - \frac{2}{T - \tau} \int_0^{T-\tau} \left(\sqrt{\frac{t}{t+\tau}}\int_0^1 \frac{2h\left(t\left(t+\tau\right)\right)^h}{h+ 1/2} \Hypergeometric{2}{1}{\frac{1}{2} - h, 1}{\frac{3}{2} + h}{\frac{t}{t+\tau}}f_\mathcal{H}\left(h\right) dh\right) dt.
   \end{aligned}
   \end{equation}
    \begin{proof}
    We start from the definition of TAMSD \eqref{tamsd}. After applying Propositions \ref{msdfbmre} and \ref{acvffbmreii} and have 
    \begin{equation}
        \begin{aligned}
            \mathbb{E}\left[\delta_\mathcal{H}\left(\tau\right)\right] = \mathbb{E}\left[\frac{1}{T - \tau} \int_0^{T-\tau} \left(B_{\mathcal{H}}^*\left(t + \tau\right) - B_{\mathcal{H}}^*\left(t\right)\right)^2 dt\right] 
            = \\ = \mathbb{E}\left[\frac{1}{T - \tau} \int_0^{T-\tau} (\left(B^*_{\mathcal{H}}\left(t + \tau\right)\right)^2 - 2B_{\mathcal{H}}^*\left(t + \tau\right)B_{\mathcal{H}}^*\left(t\right) + \left(B^*_{\mathcal{H}}\left(t\right)\right)^2 dt\right] = \\
            = \frac{1}{T - \tau} \int_0^{T-\tau} \left[M_\mathcal{H}\left(2\log \left( t + \tau\right)\right)  +  M_\mathcal{H}\left(2\log t\right)\right]dt - \\ - \frac{2}{T - \tau} \int_0^{T-\tau} \left(\sqrt{\frac{t}{t+\tau}}\int_0^1 \frac{2h\left(t\left(t+\tau\right)\right)^h}{h+ 1/2} \Hypergeometric{2}{1}{\frac{1}{2} - h, 1}{\frac{3}{2} + h}{\frac{t}{t+\tau}}f_\mathcal{H}\left(h\right) dh\right) dt.
        \end{aligned}
    \end{equation}
The change in the order of integration is made in accordance with Fubini's theorem \cite{Rudin}.
    \end{proof}
\end{proposition}
\noindent As for FBMRE, we also investigate properties of the increment process  of RL FBMRE $\{  b_{\mathcal{H}}^{*\tau}\left(t\right)\}_{t\geq 0}$, which is defined as 
\begin{equation}
    b_{\mathcal{H}}^{*\tau}\left(t\right) = B_{\mathcal{H}}^*\left(t + \tau\right) - B_{\mathcal{H}}^*\left(t\right).
\end{equation}
\noindent In the next proposition we present the form of the second moment of the increment process of RL FBMRE, which is different from the same characteristics for FBMRE (in a general case). 
\begin{proposition}
    The second moment of the increment process of RL FBMRE is given by
    \begin{equation}\label{eq1}
        \mathbb{E}\left[\left(b_{\mathcal{H}}^{*\tau}\left(t\right)\right)^2 \right] = \int_0^1 \left[\left(t+\tau\right)^{2h} + t^{2h} - \frac{4h \left(t+\tau\right)^{h - \frac{1}{2}}t^{h + \frac{1}{2}}}{h + 1/2} \Hypergeometric{2}{1}{\frac{1}{2} - h, 1}{\frac{3}{2} + h}{\frac{t}{t + \tau}}\right]f_{\mathcal{H}}\left(h\right) dh.
    \end{equation}
    \begin{proof}
To obtain \eqref{eq1} we apply the law of total expectation and Proposition \ref{incrementscovfbmiialter}. Thus, we have 
        \begin{equation}
            \begin{aligned}
                \mathbb{E}\left[\left(b_{\mathcal{H}}^{*\tau}\left(t\right)\right)^2 \right] = \mathbb{E}\left[\mathbb{E}\left[\giventhat{\left(b_{\mathcal{H}}^{*\tau}\left(t\right)\right)^2}{\mathcal{H}} \right]\right] = \\ = \int_0^1 \left[\left(t+\tau\right)^{2h} + t^{2h} - \frac{4h \left(t+\tau\right)^{h - \frac{1}{2}}t^{h + \frac{1}{2}}}{h + 1/2} \Hypergeometric{2}{1}{\frac{1}{2} - h, 1}{\frac{3}{2} + h}{\frac{t}{t + \tau}}\right] f_{\mathcal{H}}\left(h\right) dh.
            \end{aligned}
        \end{equation}
    \end{proof}
\end{proposition}
\begin{remark}
    Because the sequence of increments of RL FBMRE is not even stationary, it is obviously not ergodic.  
\end{remark}
\noindent The last considered property is the  H{\"o}lder continuity of the trajectories of RL FBMRE. 
\begin{theorem}
Let $\left\{B_{\mathcal{H}}^*\left(t\right)\right\}_{t\in\left[0,1\right]}$ be the RL FBMRE with $\mathcal{H}$ taking values from the interval $\left(H_1, 1\right)$, $0<H_1<1$. Then, there exists a version of RL FBMRE such that its trajectories are almost sure H{\"o}lder continuous for all exponents $\kappa <H_1$.
\begin{proof}
    Because absolute moments of order $q$ of FBMRE and RL FBMRE are the same, the proof is the same as in Theorem \ref{fbmrecontinuity}.
\end{proof}
\end{theorem}
\section{Two-point distribution of the Hurst exponent}\label{sec:two_point}
In this section, we discuss the properties of the FBMRE and RL FBMRE under the assumption that the Hurst exponent follows a mixture of two point distributions. For brevity in what follows we use the name "two-point distribution", which of course should not be mistakenly thought to be a transition probability. We analyze the asymptotics of the discussed characteristics for both short and long times. The main focus is on the second moments of the considered processes, as well as the second moments of their increment processes.
\\\noindent The  two-point distributed random variable $\mathcal{H}$ concentrated at two points $H_1, H_2 \in \left(0, 1\right)$, $H_1<H_2$ has the pdf given by
    \begin{equation}\label{tppdf}
        f_{\mathcal{H}}(h) = p\delta\left(h - H_1\right) + \left(1-p\right)\delta\left(h-H_2\right),
    \end{equation}
where the weight parameter $p\in \left(0,1\right)$ and $\delta\left(\cdot\right)$ is the Dirac delta.
The corresponding  moment generating function  is given by
\begin{equation}
    M_\mathcal{H}(s) = p e^{H_1 s} + \left(1-p\right)e^{H_2 s} .
\end{equation}
\subsection{Second moment of the process} 
The second moments of FBMRE and RL FBMRE in the considered case are
\begin{equation}
    \mathbb{E}\left[\left(B_{\mathcal{H}}\left(t\right)\right)^2\right] = \mathbb{E}\left[\left(B_{\mathcal{H}}^*\left(t\right)\right)^2\right] = pt^{2H_1} + \left(1-p\right)t^{2H_2}.
\end{equation}
It is easy to show that the asymptotics of the above function for short and long times take the form
\begin{equation}\label{msdfbmretpasympt}
    \mathbb{E}\left[\left(B_{\mathcal{H}}\left(t\right)\right)^2\right] = \mathbb{E}\left[\left(B_{\mathcal{H}}^*\left(t\right)\right)^2\right] \sim
    \begin{cases}
        pt^{2H_1}, ~~\text{ for } ~~t<<1, \\ 
        \left(1-p\right)t^{2H_2}, ~~\text{ for } ~~t>>1.
    \end{cases}
\end{equation}

\noindent The effect of changing exponent in asymptotics of the second moment is called the accelerating diffusion \cite{Chechkin2002, Chechkin2008, Chechkin2011}.
\subsection{Autocovariance function of the process}
\noindent In the considered case the autocovariance function of FBMRE has the following form
\begin{equation}\label{covtp}
    Cov\left(B_{\mathcal{H}}\left(t\right), B_{\mathcal{H}}\left(t+\tau\right)\right) = \frac{1}{2}\left[p\left((t+\tau)^{2H_1} + t^{2H_1} - \tau^{2H_1}\right) + \left(1-p\right)\left((t+\tau)^{2H_2} + t^{2H_2} - \tau^{2H_2}\right) \right],
\end{equation}
while for RL FBMRE we have
\begin{equation}\label{coviitp}
    \begin{aligned}
    Cov\left(B_{\mathcal{H}}^*\left(t\right), B_{\mathcal{H}}^*\left(t+\tau\right)\right) = \sqrt{\frac{t}{t+\tau}}\left[p\left(t\left(t+\tau\right)\right)^{H_1} \Hypergeometric{2}{1}{\frac{1}{2} - H_1, 1}{\frac{3}{2} + H_1}{\frac{t}{t+\tau}}\right. + \\ +\left.\left(1-p\right)\left(t\left(t+\tau\right)\right)^{H_2} \Hypergeometric{2}{1}{\frac{1}{2} - H_2, 1}{\frac{3}{2} + H_2}{\frac{t}{t+\tau}}\right].
    \end{aligned}
\end{equation}
Applying Taylor expansion we obtain the following asymptotics of the autocovariance function of FBMRE
\begin{equation}\label{covfbmretpasympt}
    Cov\left(B_{\mathcal{H}}\left(t\right), B_{\mathcal{H}}\left(t+\tau\right)\right) \sim
    \begin{cases}
        \frac{1}{2}\left[p\tau^{2H_1}\left(2H_1\frac{t}{\tau}  + \left(\frac{t}{\tau}\right)^{2H_1} + H_1\left(2H_1-1\right)\left(\frac{t}{\tau}\right)^2 \right) + \right. \\ \left. + \left(1-p\right)\tau^{2H_2}\left(2H_2\frac{t}{\tau}  + \left(\frac{t}{\tau}\right)^{2H_2} + H_2\left(2H_2-1\right)\left(\frac{t}{\tau}\right)^2\right)\right], ~~\text{ for } ~~t/\tau<<1\\
        \frac{1}{2}\left[pt^{2H_1}\left(2 + 2H_1\frac{\tau}{t}  - \left(\frac{\tau}{t}\right)^{2H_1} + H_1\left(2H_1-1\right)\left(\frac{\tau}{t}\right)^2\right) + \right. \\ \left. + \left(1-p\right)t^{2H_2}\left(2 + 2H_2\frac{\tau}{t} + H_2\left(2H_2-1\right)\left(\frac{\tau}{t}\right)^2 - \left(\frac{\tau}{t}\right)^{2H_2}\right)\right], ~~\text{ for }~~ t/\tau>>1.
    \end{cases}
\end{equation}
Using Proposition \ref{fbmiicovasympt} one can show that  for  $t/\tau << 1$, we have the following asymptotics of the autocovariance function of RL FBMRE
\begin{equation}\label{covfbmreiitpasymptshort}
    \begin{aligned}
    Cov\left(B_{\mathcal{H}}^*\left(t\right), B_{\mathcal{H}}^*\left(t+\tau\right)\right) \sim
        p\tau^{2H_1}\left(\frac{t}{\tau}\right)^{H_1+\frac{1}{2}} \left[\frac{2H_1}{\frac{1}{2} + H_1} - \frac{2H_1\left(\frac{1}{2}-H_1\right)}{\frac{3}{2} + H_1}\frac{t}{\tau} +  \frac{H_1\left(\frac{1}{2}-H_1\right)\left(\frac{3}{2}-H_1\right)}{\frac{5}{2} + H_1}\left(\frac{t}{\tau}\right)^{2}\right] + \\ + \left(1-p\right)\tau^{2H_2}\left(\frac{t}{\tau}\right)^{H_2+\frac{1}{2}} \left[\frac{2H_2}{\frac{1}{2} + H_2} - \frac{2H_2\left(\frac{1}{2}-H_2\right)}{\frac{3}{2} + H_2}\frac{t}{\tau} +  \frac{H_2\left(\frac{1}{2}-H_2\right)\left(\frac{3}{2}-H_2\right)}{\frac{5}{2} + H_2}\left(\frac{t}{\tau}\right)^{2}\right].
    \end{aligned}
\end{equation}
while for $t/\tau>>1$ the following holds
\begin{equation}\label{covfbmreiitpasymptlong}
    \begin{aligned}
        Cov\left(B_{\mathcal{H}}^*\left(t\right), B_{\mathcal{H}}^*\left(t+\tau\right)\right) \sim
         \frac{1}{2}\left[pt^{2H_1}\left(2 + 2H_1\frac{\tau}{t}  - \frac{2H_1\Gamma\left(H_1 + \frac{1}{2}\right)^2}{\Gamma\left(1+2H_1\right)\sin\left(\pi H_1\right)}\left(\frac{\tau}{t}\right)^{2H_1} + H_1\left(2H_1-1\right)\left(\frac{\tau}{t}\right)^2 \right) + \right. \\ \left. + \left(1-p\right)t^{2H_2}\left(2 + 2H_2\frac{\tau}{t}  - \frac{2H_2\Gamma\left(H_2 + \frac{1}{2}\right)^2}{\Gamma\left(1+2H_2\right)\sin\left(\pi H_2\right)}\left(\frac{\tau}{t}\right)^{2H_2} + H_2\left(2H_2-1\right)\left(\frac{\tau}{t}\right)^2\right) \right].
    \end{aligned}
\end{equation}
Comparing the behavior of the autocovariance of FBMRE and RL FBMRE for $t/\tau<<1$, we observe that the orders of the dominant terms in both cases are different. Specifically, for FBMRE, the dominant terms are of order $\mathcal{O}\left(\tau^{2H_1}\frac{t}{\tau}\right)$, whereas for RL FBMRE, the dominant terms are of order $\mathcal{O}\left(\tau^{2H_1}\left(\frac{t}{\tau}\right)^{H_1+\frac{1}{2}}\right)$. However, for $t/\tau>>1$, we find that in both cases, the dominant terms are of order  $\mathcal{O}\left(t^{2H_1}\right)$. The situation is similar to what we have in case of constant Hurst exponent in equations \eqref{fbmiicovasympteq} and \eqref{covfbmtpasympt}.

\begin{figure}[ht!]
       \centering
         \includegraphics[width=1\textwidth, height=0.3\textheight]{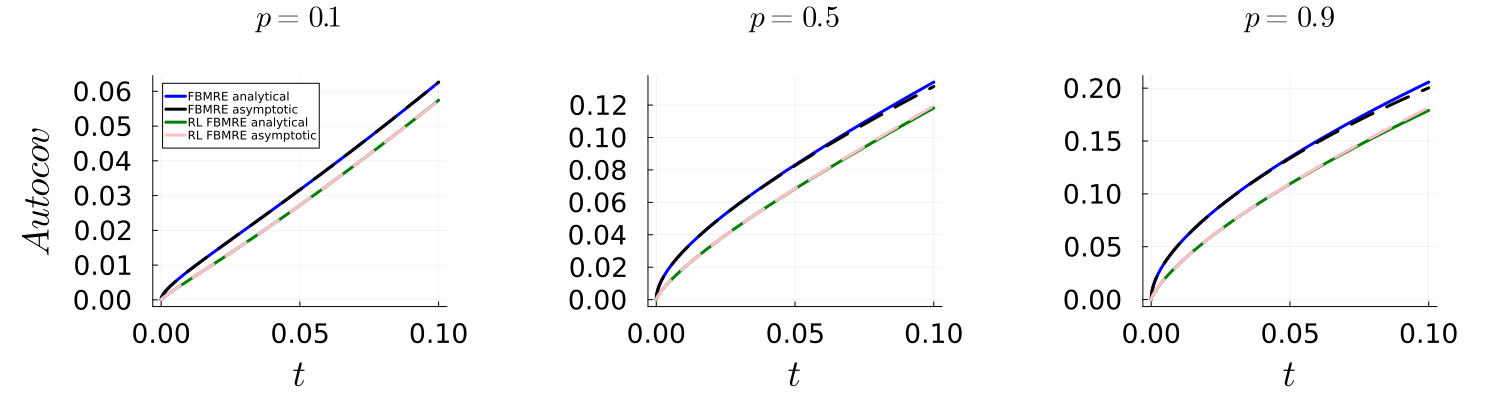}
        \caption{Comparison of the autocovariance functions obtained from analytical formulas and their asymptotic behavior for $t/\tau<<1$ for FBMRE and RL FBMRE, assuming  two-point distribution of the Hurst exponent with $H_1=0.25$, $H_2=0.75$. Left panel:  $p=0.1$; middle panel: $p=0.5$; right panel $p=0.9$. We assume  $\tau=0.1$.}\label{fig::covshort}
\end{figure}

\begin{figure}[ht!]
       \centering
         \includegraphics[width=1\textwidth, height=0.3\textheight]{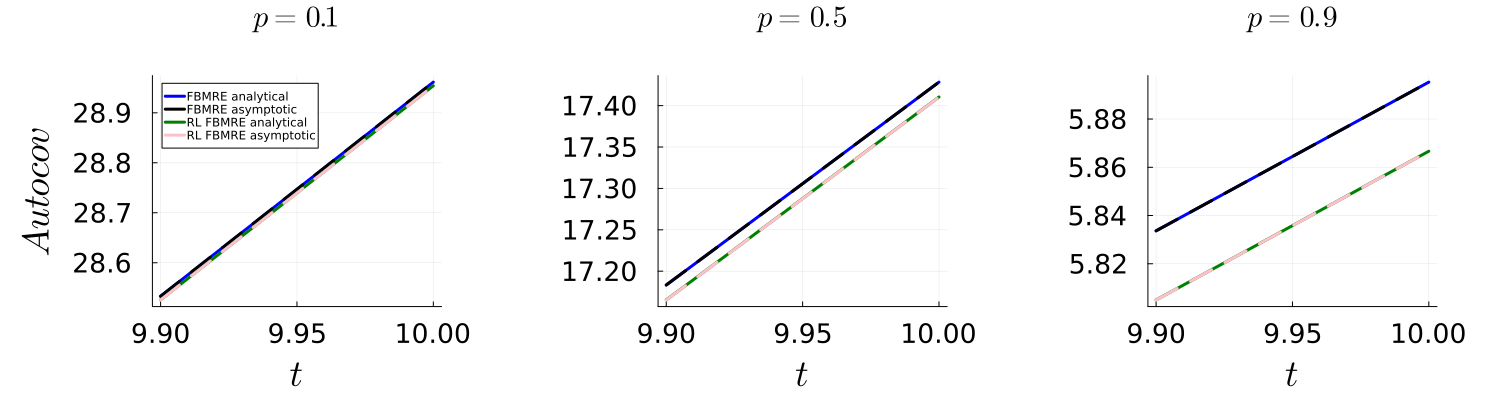}
        \caption{Comparison of the autocovariance functions obtained from analytical formulas and their asymptotic behavior for $t/\tau>>1$ for FBMRE and RL FBMRE, assuming  two-point distribution of the Hurst exponent with $H_1=0.25$, $H_2=0.75$. Left panel:  $p=0.1$; middle panel: $p=0.5$; right panel $p=0.9$. We assume  $\tau=0.1$.}\label{fig::covlong}
\end{figure}

To illustrate the theoretical results discussed above, in Figs. \ref{fig::covshort} and \ref{fig::covlong} we present the autocovariance function for FBMRE and RL FBMRE for short and long times, respectively. On the plots we present the curves corresponding to equations in the text. We assume that the Hurst exponent follows a two-point distribution with $H_1=0.25$ and $H_2=0.75$. The figures include separate subplots demonstrating the behavior of the autocovariance functions for $p=0.1$ (left panels), $p=0.5$ (middle panels) and $p=0.9$ (right panels). Additionally, we assume $\tau=0.1$. The autocovariance functions obtained from the analytical formulas (\ref{covtp}) and (\ref{coviitp}) for FBMRE and RL FBMRE, respectively, are presented along with their asymptotics (\ref{covfbmretpasympt}) for FBMRE, and (\ref{covfbmreiitpasymptshort}) for RL FBMRE for short times, and (\ref{covfbmreiitpasymptlong}) for RL FBMRE for long times.

A perfect agreement is observed between the analytical formulas and the corresponding derived asymptotics. However, it is important to emphasize that this agreement is achieved only if all terms contained in the corresponding asymptotic formulas are included.  The difference between the autocovariance functions of the two processes is clearly visible for all $p$ values considered at both short and long times. However, as mentioned earlier, the dominant factor of the autocovariance functions of both processes at long times is the same.

\subsection{Expectation of TAMSD}
\noindent For the FBMRE we have the following formula for the expectation of the TAMSD
\begin{equation}\label{emtamsdtp}
    \mathbb{E}\left[\delta_\mathcal{H}\left(\tau\right)\right] = p\tau^{2H_1} + \left(1-p\right)\tau^{2H_2},
\end{equation}
while for RL FBMRE we have
\begin{equation}\label{etamsdiitp}
    \begin{aligned}
   \mathbb{E}\left[\delta_\mathcal{H}^*\left(\tau\right)\right] = p\left(\frac{1}{2H_1+1}\frac{T^{2H_1+1} - \tau^{2H_1+1}}{T- \tau} + \frac{1}{2H_1+1}\left(T - \tau\right)^{2H_1} -\right. \\ \left. - \frac{4H_1}{\left(H_1+\frac{1}{2}\right)\left(H_1+\frac{3}{2}\right)} \frac{\left(T - \tau\right)^{H_1 + 1/2}}{\tau^{1/2 - H_1}}\Hypergeometric{2}{1}{\frac{1}{2} + H_1, \frac{1}{2}-H_1}{\frac{5}{2} + H_1}{\frac{\tau-T}{\tau}}\right) + \\ + \left(1-p\right)\left(\frac{1}{2H_2+1}\frac{T^{2H_2+1} - \tau^{2H_2+1}}{T- \tau} + \frac{1}{2H_2+1}\left(T - \tau\right)^{2H_2} - \right. \\ \left. - \frac{4H_2}{\left(H_2+\frac{1}{2}\right)\left(H_2+\frac{3}{2}\right)} \frac{\left(T - \tau\right)^{H_2 + 1/2}}{\tau^{1/2 - H_2}}\Hypergeometric{2}{1}{\frac{1}{2} + H_2, \frac{1}{2}-H_2}{\frac{5}{2} + H_2}{\frac{\tau-T}{\tau}}\right).
    \end{aligned}
\end{equation}
One can show that in case of FBMRE and for $T/\tau>>1$ the following holds
\begin{equation}\label{tamsdfbmretpasympt}
    \mathbb{E}\left[\delta_\mathcal{H}\left(\tau\right)\right] \sim
    \begin{cases}
        p\tau^{2H_1}, ~~\text{ for }~~ \tau<<1, \\ 
        \left(1-p\right)\tau^{2H_2}, ~~\text{ for } ~~\tau>>1.
    \end{cases}
\end{equation}
For RL FBMRE and in case  $T/\tau>>1$, by the use of Theorem \ref{tamsdfbmiiasympt} we obtain 
\begin{equation}
        \begin{aligned}
            \mathbb{E}\left[\delta_{\mathcal{H}}^*\left(\tau\right)\right] \sim \begin{cases}

                p\frac{2H_1\Gamma\left(\frac{1}{2}+H_1\right)^2}{\Gamma\left(1+2H_1\right)\sin\left(\pi H_1\right)}\tau^{2H_1}, ~~\text{ for } ~~\tau<<1, \\ 
                \left(1-p\right)\frac{2H_2\Gamma\left(\frac{1}{2}+H_2\right)^2}{\Gamma\left(1+2H_2\right)\sin\left(\pi H_2\right)} \tau^{2H_2}, ~~\text{ for }\frac{1}{2}, \, \tau>>1.
    \end{cases}
        \end{aligned}
    \end{equation}
Despite the differences in the exact formulas for FBMRE and RL FBMRE, under the assumption that $T/\tau>>1$, we observe similar asymptotic behavior in both cases. Specifically, up to a constant, for $\tau<<1$, the expectation of TAMSD behaves like $\tau^{2H_1}$, and for $\tau>>1$, it behaves like $\tau^{2H_2}$.

\begin{figure}[h!]
       \centering
         \includegraphics[width=1\textwidth, height=0.3\textheight]{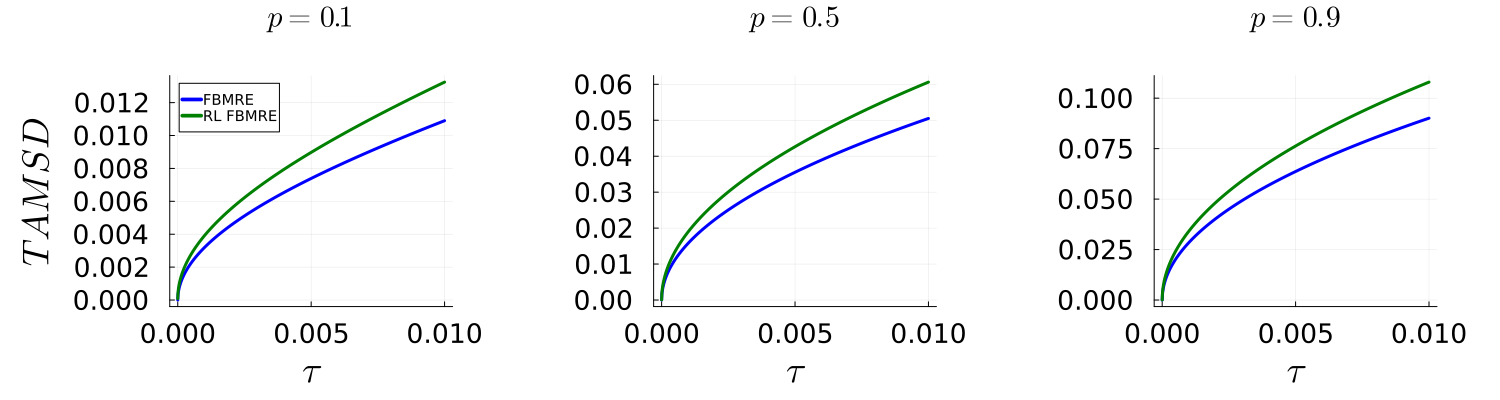}
        \caption{Comparison of the expectations of TAMSD for   FBMRE and RL FBMRE obtained from the analytical formulas, assuming  two-point distribution of the Hurst exponent with $H_1=0.25$ and $H_2=0.75$. Left panel:  $p=0.1$; middle panel: $p=0.5$; right panel $p=0.9$. We assume  $T/\tau>>1$, $\tau<<1$ and $T=20000$.}\label{fig::eatamsdshort}
\end{figure}

\begin{figure}[h!]
       \centering
         \includegraphics[width=1\textwidth, height=0.3\textheight]{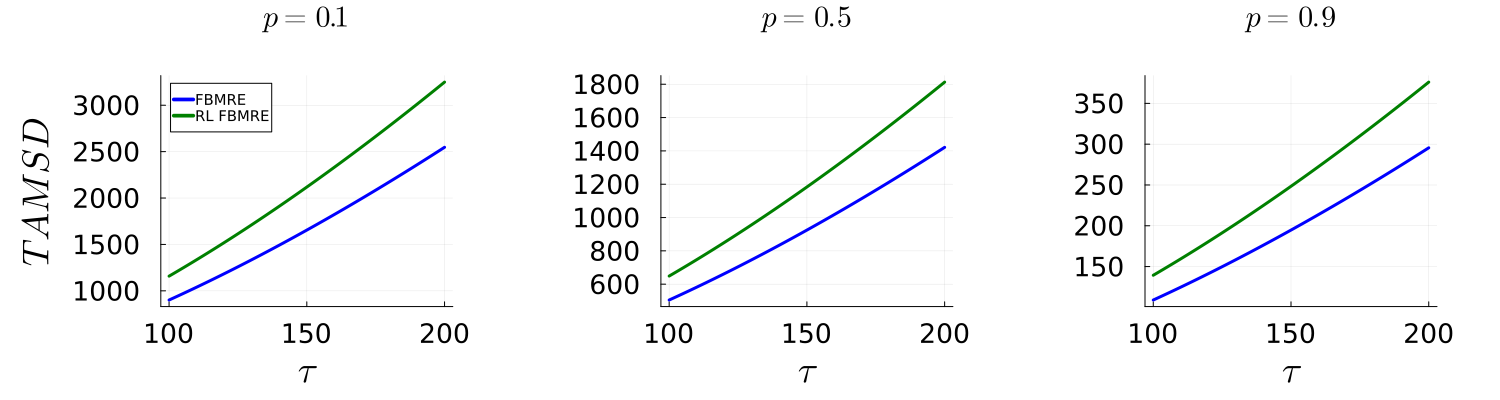}
        \caption{Comparison of the expectations of TAMSD for   FBMRE and RL FBMRE obtained from the analytical formulas, assuming  two-point distribution of the Hurst exponent with $H_1=0.25$ and $H_2=0.75$. Left panel:  $p=0.1$; middle panel: $p=0.5$; right panel $p=0.9$. We assume  $T/\tau>>1$, $\tau>>1$ and $T=20000$.}\label{fig::eatamsdlong}
\end{figure}

To illustrate the difference between the expectations of the TAMSD for FBMRE and RL FBMRE with a two-point distributed Hurst exponent, in Figs. \ref{fig::eatamsdshort} and \ref{fig::eatamsdlong} we present the discussed characteristic  obtained from the analytical formulas (\ref{emtamsdtp}) and (\ref{etamsdiitp}) for FBMRE and RL FBMRE, respectively. Specifically, we assume $T/\tau>>1$ and $T=20000$; the case $\tau<<1$ is shown in Fig. \ref{fig::eatamsdshort}, while the case $\tau>>1$ is illustrated in Fig. \ref{fig::eatamsdlong}. As in the illustration of the autocovariance function, we assume $H_1=0.25$ and $H_2=0.75$, and in each subplot, we present the behavior of the expectation of TAMSD for different values of $p$: $p=0.1$ (left panels), $p=0.5$ (middle panels) and $p=0.9$ (right panels). In contrast to the autocovariance function, we observe a significant difference in this statistic for both processes, making it more effective in highlighting the essential differences between them.


\subsection{Second moment of the increment process}
For the two-point distribution of the Hurst exponent, the second moment of increment process of FBMRE is equal to the expectation of the TAMSD. Thus, it is given by
\begin{equation}\label{tpincrementsvarfbmre}
    \mathbb{E}\left[\left(b_{\mathcal{H}}^{\tau}\left(t\right)\right)^2 \right] = p\tau^{2H_1} + \left(1-p\right)\tau^{2H_2}.
\end{equation}
For RL FBMRE we have the following
\begin{equation}\label{tpincrementsvarfbmreii}
        \begin{aligned}
            \mathbb{E}\left[\left(b_{\mathcal{H}}^{*\tau}\left(t\right)\right)^2 \right] = p\left[\left(t+\tau\right)^{2H_1} + t^{2H_1} - \frac{4H_1 \left(t+\tau\right)^{H_1 - \frac{1}{2}}t^{H_1 + \frac{1}{2}}}{H_1 + 1/2} \Hypergeometric{2}{1}{\frac{1}{2} - H_1, 1}{\frac{3}{2} + H_1}{\frac{t}{t + \tau}}\right]  + \\ + \left(1-p\right)\left[\left(t+\tau\right)^{2H_2} + t^{2H_2} - \frac{4H_2 \left(t+\tau\right)^{H_2 - \frac{1}{2}}t^{H_2 + \frac{1}{2}}}{H_2 + 1/2} \Hypergeometric{2}{1}{\frac{1}{2} - H_2, 1}{\frac{3}{2} + H_2}{\frac{t}{t + \tau}}\right],
        \end{aligned}
\end{equation}
which gives the same asymptotic, as for the expectation of TAMSD, see Eq. (\ref{tamsdfbmretpasympt}).
\\
\noindent Using Lemma \ref{Isampt} one gets the following asymptotic for RL FBMRE in case  $t/\tau<<1$ 
\begin{equation}\label{varincfbmreiitpasymptshort}
    \mathbb{E}\left[\left(b_{\mathcal{H}}^{*\tau}\left(t\right)\right)^2 \right] \sim
         p\tau^{2H_1} + \left(1-p\right)\tau^{2H_2}, 
\end{equation}
while for $t/\tau>>1$ it reads as \begin{equation}\label{varincfbmreiitpasymptlonh}
    \mathbb{E}\left[\left(b_{\mathcal{H}}^{*\tau}\left(t\right)\right)^2 \right] \sim
         p\tau^{2H_1} \frac{2H_1\Gamma\left(H_1 + \frac{1}{2}\right)^2}{\Gamma\left(1+2H_1\right)\sin\left(\pi H_1\right)} + \left(1-p\right)\tau^{2H_2} \frac{2H_2\Gamma\left(H_2 + \frac{1}{2}\right)^2}{\Gamma\left(1+2H_2\right)\sin\left(\pi H_2\right)}.
\end{equation}
 For $t/\tau<<1$, we observe the same asymptotics of the second moments of FBMRE and RL FBMRE. Namely, the dominant terms in both cases are of order $\mathcal{O}\left(\tau^{2H_1}\right)$, when $\tau<<1$ and of order $\mathcal{O}\left(\tau^{2H_2}\right)$, when $\tau>>1$. However,when $t/\tau>>1$, the dominant terms in both cases are also of order $\mathcal{O}\left(\tau^{2H_1}\right)$, when $\tau<<1$ and of order $\mathcal{O}\left(\tau^{2H_2}\right)$, when $\tau>>1$, but the prefactors are different. This is the difference between these two processes, which can be detected in experiment.

\section{Conclusions}
In this paper, we explored the probabilistic properties of RL FBMRE, the Riemann-Liouville fractional Brownian motion with a random Hurst exponent. This process retains the features of RL FBM at the individual trajectory level, but the Hurst exponent varies randomly across the trajectories. RL FBM, like FBM, is a Gaussian process that exhibits anomalous diffusion, though its increments are non-stationary. Besides their Gaussian nature, RL FBM and FBM have the same second moment. However, they differ in certain probabilistic properties, including the expectation of the TAMSD and the autocovariance function. Our study was motivated by recent biophysical experiments showing that the anomalous diffusion models with fixed parameters are inadequate for describing the observed phenomena. This naturally requires modification of FBM and RL FBM by introducing random parameters, a scenario previously discussed for FBM. While the case of FBM with random parameters has been studied in the literature, to the best of our knowledge, the probabilistic properties of RL FBMRE were not studied. Firstly, we made significant advances in the theory of FBMRE by calculating its autocovariance, expectation of TAMSD and the second moment of increments. Moreover, we proved non-ergodicity of the sequence of increments of FBMRE, in contrast to the increments of FBM. Remarkably, despite the lack of ergodicity the second moment of FBMRE and the expectation of its TAMSD are the same. We also analyzed the  H{\"o}lder continuity of trajectories of FBMRE. Secondly, we established the mathematical framework for RL FBMRE by analyzing its second moment,  autocovariance function, ergodicity property, expectation of TAMSD, H{\"o}lder continuity of its trajectories and second moment of the incremental process. As a specific case of the anomalous diffusion exponent distribution, we considered the mixture of two point distributions. For this case, we analyzed the asymptotic behavior of the characteristics of FBMRE and RL FBMRE. Our theoretical findings have been confirmed through numerical analysis. We observed that the autocovariance functions for FBMRE and RL FBMRE exhibit different behavior for small time-to-lag ratios, while they are the same in the opposite limit. However, for the expectations of TAMSD we observed different behavior for both short and long times. In case of the second moment of the increments, similarly to the autocovariance functions, the behavior of both processes differs for small time-to-lag ratios but is the same otherwise.

The theoretical results presented elucidate the principal differences between the two examined processes with a random Hurst exponent. These findings also provide a foundation for proposing new techniques to distinguish such processes, which is analogous to the approaches developed in \cite{doi:10.1063/5.0044878} for Gaussian processes and also in \cite{10.1063/5.0201436}, where the authors introduced a statistical methodology for differentiating FBM from FBMRE based on time-averaged statistics.

\section*{Acknowledgements}
We thank Qing Wei (Chinese Academy of Science) and Wei Wang (University of Potsdam) for helpful discussions and remarks about derivation of expectation of TAMSD and its asymptotics.
Agnieszka Wy\l oma\'nska acknowledges support from the National Science Centre, Poland, via project No. 2020/37/B/HS4/00120. Aleksei Chechkin acknowledges support from the BMBF Project PLASMA-SPIN-ENERGY.
\bibliography{bibliography.bib}
\appendix
\section{Riemann-Liouville ractional Brownian motion - proofs}\label{AppA}
In the proposition presented below we recall a well known fact that will be usefu in the further analysis.
\begin{proposition}\label{Isampt}
    \noindent The function $I\left(t, \tau; H\right)$ defined in Eq. (\ref{II}) has the following  asymptotic for  $t/\tau>>1$\cite{Balcerek2022}
       \begin{equation}
I\left(t, \tau; H\right) + \frac{1}{2H} \sim
         \frac{\Gamma\left(H + \frac{1}{2}\right)^2\Gamma\left(1-2H\right)\cos\left(\pi H\right)}{ \pi H} = \frac{\Gamma\left(H + \frac{1}{2}\right)^2}{\Gamma\left(1+2H\right)\sin\left(\pi H\right)}.
\end{equation}
\end{proposition}

\noindent \textbf{Proof of Theorem \ref{th1} }
\begin{proof} We have the following
    \begin{equation}\label{calculatingtamsd}
    \begin{aligned}
\mathbb{E}\left[\delta_H^*\left(\tau\right)\right] = \mathbb{E}\left[\frac{1}{T - \tau} \int_0^{T-\tau} \left(B_H^*\left(t + \tau\right) - B_H^*\left(t\right)\right)^2 dt\right]  
    =\\=  \frac{1}{T - \tau} \int_0^{T-\tau} \mathbb{E}\left[\left(B^*_H\left(t + \tau\right)\right)^2 - 2B_H^*\left(t + \tau\right)B_H^*\left(t\right) + \left(B^*_H\left(t\right)\right)^2\right]  dt = \\ = \frac{1}{T - \tau} \int_0^{T-\tau} \left[\left(t + \tau\right)^{2H} + t^{2H} - \frac{4H}{H+\frac{1}{2}}\sqrt{\frac{t}{t+\tau}}\left(\left(t\left(t+\tau\right)\right)^{H} \Hypergeometric{2}{1}{\frac{1}{2} - H, 1}{\frac{3}{2} + H}{\frac{t}{t+\tau}}\right)\right] dt = \\ = 
    \frac{1}{2H+1}\frac{T^{2H+1} - \tau^{2H+1}}{T- \tau} + \frac{1}{2H+1}\left(T - \tau\right)^{2H} - \\ - \frac{4H}{H+\frac{1}{2}} \frac{1}{T-\tau} \int_0^{T-\tau} t^{H + 1/2}\left(t+\tau\right)^{H - 1/2}\Hypergeometric{2}{1}{\frac{1}{2} - H, 1}{\frac{3}{2} + H}{\frac{t}{t+\tau}}dt = \\ = \frac{1}{2H+1}\frac{T^{2H+1} - \tau^{2H+1}}{T- \tau} + \frac{1}{2H+1}\left(T - \tau\right)^{2H} - I_H\left(\tau, T\right).
    \end{aligned}
    \end{equation}
    Change of order of integration and expectation is performed according to the Fubinis theorem. Now, we will calculate $I_H\left(\tau, T\right)$. Using formula (15.3.4) from \cite{Abramowitz}, the relationship between hypergeometric functions and Fox H-function given by formula (1.131) from \cite{Mathai2010}, then applying formulas (1.16.4.1), (8.3.2.1), (8.3.2.6)  and (8.3.2.7) from \cite{Prudnikov} we have
    \begin{equation}\label{I1}
    \begin{aligned}
         I_H\left(\tau, T\right) =  \frac{1}{T-\tau}\frac{4H}{H+\frac{1}{2}}  \int_0^{T-\tau} t^{H + 1/2}\tau^{H - 1/2}\Hypergeometric{2}{1}{\frac{1}{2} - H, H + \frac{1}{2}}{\frac{3}{2} + H}{-\frac{t}{\tau}}dt = \\ = \frac{4H}{\Gamma\left(\frac{1}{2}-H\right)} \frac{\tau^{2H + 1}}{T - \tau} \int_0^{\frac{T-\tau}{\tau}} s^{H+1/2}\FoxH{1}{2}{2}{2}{s}{H + \frac{1}{2}, \frac{1}{2} - H}{1, 1}{0, -H - \frac{1}{2}}{1, 1}ds = \\ = \frac{4H}{\Gamma\left(\frac{1}{2}-H\right)} \frac{\left(T - \tau\right)^{H + 1/2}}{\tau^{1/2 - H}} \FoxH{1}{2}{2}{2}{\frac{T-\tau}{\tau}}{H+\frac{1}{2}, \frac{1}{2} - H}{1,0}{0, -H - \frac{3}{2}}{1, 1} = \\ = \frac{4H}{\left(H+\frac{1}{2}\right)\left(H+\frac{3}{2}\right)} \frac{\left(T - \tau\right)^{H + 1/2}}{\tau^{1/2 - H}}\Hypergeometric{2}{1}{\frac{1}{2} + H, \frac{1}{2}-H}{\frac{5}{2} + H}{\frac{\tau-T}{\tau}}.
    \end{aligned}
    \end{equation}
Putting \eqref{I1} into \eqref{calculatingtamsd} we obtain \eqref{fbmiitamsd}.
\end{proof}
\noindent \textbf{Proof of Theorem \ref{tamsdfbmiiasympt}}
\begin{proof}
    We need to approximate the three terms in Eq. \eqref{fbmiitamsd}. Using (8.3.2.7) from \cite{Prudnikov} we have 
    \begin{equation}
    \begin{aligned}
        \mathbb{E}\left[\delta_H^*\left(\tau\right)\right] =  
            \frac{1}{2H+1}\frac{T^{2H+1} - \tau^{2H+1}}{T- \tau} + \frac{1}{2H+1}\left(T - \tau\right)^{2H} - \\ - \frac{4H}{\Gamma\left(\frac{1}{2}-H\right)} \frac{\left(T - \tau\right)^{H + 1/2}}{\tau^{1/2 - H}} \FoxH{1}{2}{2}{2}{\frac{T-\tau}{\tau}}{H+\frac{1}{2}, \frac{1}{2} - H}{1,0}{0, -H - \frac{3}{2}}{1, 1} = \\ = 
            \frac{1}{2H+1}\frac{T^{2H+1} - \tau^{2H+1}}{T- \tau} + \frac{1}{2H+1}\left(T - \tau\right)^{2H} - \\ - \frac{4H}{\Gamma\left(\frac{1}{2}-H\right)} \frac{\left(T - \tau\right)^{H + 1/2}}{\tau^{1/2 - H}} \FoxH{2}{1}{2}{2}{\frac{\tau}{T-\tau}}{1, H + \frac{5}{2}}{1,1}{  \frac{1}{2} - H, H + \frac{1}{2}}{1, 1} = I_1 + I_2 - I_{H}\left(\tau, T\right).
    \end{aligned}
    \end{equation}
First, we approximate $I_1$ and $I_2$ with the help of the Taylor expansion considering only its first three components.
Finally, we approximate $ I_{H}\left(\tau, T\right)$. Using formula (8.3.2.3) from \cite{Prudnikov} we obtain the series representation for $H\neq\frac{1}{2}$,
\begin{equation}
    \begin{aligned}
         I_{H}\left(\tau, T\right) = 
         \frac{4H}{\Gamma\left(\frac{1}{2}-H\right)} \frac{\left(T - \tau\right)^{H + 1/2}}{\tau^{1/2 - H}}\sum_{k=0}^{\infty}\Gamma\left(2H-k\right)\frac{\Gamma\left(\frac{1}{2}-H+k\right)}{\Gamma\left(2H+2-k\right)}\frac{\left(-1\right)^k \left(\frac{\tau}{T-\tau}\right)^{\frac{1}{2}-H+k}}{k!} + \\ + \frac{4H}{\Gamma\left(\frac{1}{2}-H\right)} \frac{\left(T - \tau\right)^{H + 1/2}}{\tau^{1/2 - H}}\sum_{k=0}^{\infty}\Gamma\left(-2H-k\right)\frac{\Gamma\left(\frac{1}{2}+H+k\right)}{\Gamma\left(2-k\right)}\frac{\left(-1\right)^k \left(\frac{\tau}{T-\tau}\right)^{\frac{1}{2}+H+k}}{k!} = \\ 
         = \frac{4H T^{2H}}{\Gamma\left(\frac{1}{2}-H\right)} \sum_{k=0}^{\infty}\frac{\Gamma\left(\frac{1}{2}-H+k\right)\left(-1\right)^k }{\left(2H-k\right)\left(2H+1-k\right)k!}\left(\frac{\tau}{T}\right)^k\left(1-\frac{\tau}{T}\right)^{2H-k} + \\ + \frac{4H\tau^{2H}}{\Gamma\left(\frac{1}{2}-H\right)}\sum_{k=0}^{\infty}\Gamma\left(-2H-k\right)\frac{\Gamma\left(\frac{1}{2}+H+k\right)\left(-1\right)^k }{\Gamma\left(2-k\right)k!}\left(\frac{\tau}{T-\tau}\right)^{k}.
    \end{aligned}
\end{equation}
In the first series, we expand $\left(1-\frac{\tau}{T}\right)^{2H-k}$ using Taylor series. Let us note, that in the second series of the above formula for $k\geq 2$ the gamma function diverges to infinity, and thus the appropriate terms vanish.  Therefore we have 
\begin{equation}\label{ihest}
\begin{aligned}
           I_{H}\left(\tau, T\right) 
            \sim \frac{4H T^{2H}}{\Gamma\left(\frac{1}{2}-H\right)} \sum_{k=0}^{\infty}\frac{\Gamma\left(\frac{1}{2}-H+k\right)\left(-1\right)^k }{\left(2H-k\right)\left(2H+1-k\right)k!}\left(\frac{\tau}{T}\right)^k\left[1 - \left(2H-k\right)\frac{\tau}{T} + \frac{1}{2}\left(2H-k\right)\left(2H-k-1\right)\left(\frac{\tau}{T}\right)^2\right] +\\ + \frac{4H\tau^{2H}}{\Gamma\left(\frac{1}{2}-H\right)}\left[\Gamma\left(-2H\right)\Gamma\left(\frac{1}{2}+H\right) - \Gamma\left(-2H-1\right)\Gamma\left(\frac{3}{2}+H\right)\left(\frac{\tau}{T-\tau}\right)\right] 
            \sim \\ \sim \frac{4H T^{2H}}{\Gamma\left(\frac{1}{2}-H\right)} \left[\frac{\Gamma\left(\frac{1}{2}-H\right)}{2H\left(2H+1\right)}\left(1 - 2H\frac{\tau}{T} + \frac{1}{2}2H\left(2H-1\right)\left(\frac{\tau}{T}\right)^2\right) - \right. \\  -\frac{\Gamma\left(\frac{3}{2}-H\right) }{\left(2H-1\right)2H}\frac{\tau}{T}\left(1 - \left(2H-1\right)\frac{\tau}{T} + \frac{1}{2}\left(2H-1\right)\left(2H-2\right)\left(\frac{\tau}{T}\right)^2\right) + \\   \left. + \frac{\Gamma\left(\frac{5}{2}-H\right) }{2\left(2H-2\right)\left(2H-1\right)}\left(\frac{\tau}{T}\right)^2\left(1 - \left(2H-2\right)\frac{\tau}{T} + \frac{1}{2}\left(2H-2\right)\left(2H-3\right)\left(\frac{\tau}{T}\right)^2\right)\right] + \\ + \frac{4H\tau^{2H}\Gamma\left(-2H\right)\Gamma\left(\frac{1}{2}+H\right)}{\Gamma\left(\frac{1}{2}-H\right)}\left[1 + \frac{1}{2}\left(\frac{\tau}{T-\tau}\right)\right] 
            \sim \\ \sim 2H T^{2H} \left[\frac{1}{H\left(2H+1\right)} + \frac{\frac{1}{2}-H}{H\left(2H+1\right)}\frac{\tau}{T}+ \left(\frac{\frac{1}{2}-H}{H\left(2H+1\right)} + \frac{H - \frac{3}{2}}{2\left(2H-2\right)}\right)\left(\frac{\tau}{T}\right)^2\right] + \\ + \frac{4H\Gamma\left(-2H\right)\Gamma\left(\frac{1}{2}+H\right)T^{2H}}{\Gamma\left(\frac{1}{2}-H\right)}\left[\left(\frac{\tau}{T}\right)^{2H} + \frac{1}{2}\left(\frac{\tau}{T}\right)^{2H+1}\right].
\end{aligned}
\end{equation}
Observing that the main terms in Talor expansions of $I_1$ and $I_2$ cancel out with main terms in \eqref{ihest}, we obtain the following
\begin{equation}\label{etamsdfbmiiest1}
\begin{aligned}
\mathbb{E}\left[\delta_H^*\left(\tau\right)\right] \sim  \frac{H\left(H - \frac{1}{2}\right)T^{2H}}{2H-2}\left(\frac{\tau}{T}\right)^2 + \frac{4H\pi\Gamma\left(\frac{1}{2}+H\right)T^{2H}}{\Gamma\left(\frac{1}{2}-H\right)\Gamma\left(1+2H\right)\sin\left(2\pi H\right)}\left[\left(\frac{\tau}{T}\right)^{2H} + \frac{1}{2}\left(\frac{\tau}{T}\right)^{2H+1}\right] = \\ = \frac{H\left(H - \frac{1}{2}\right)T^{2H}}{2H-2}\left(\frac{\tau}{T}\right)^2 + \frac{2H\Gamma\left(\frac{1}{2}+H\right)^2T^{2H}}{\Gamma\left(1+2H\right)\sin\left(\pi H\right)}\left[\left(\frac{\tau}{T}\right)^{2H} + \frac{1}{2}\left(\frac{\tau}{T}\right)^{2H+1}\right] .
\end{aligned}
\end{equation}
If $H=\frac{1}{2}$, then RL FBM reduces to standard Brownian motion and $\mathbb{E}\left[\delta_H^*\left(\tau\right)\right] = \tau$. Finally, keeping only the dominating terms in \eqref{etamsdfbmiiest1} we obtain \eqref{tamsdfbmiiasympteq}.
\end{proof}
\noindent \textbf{Proof of Proposition \ref{incvarfbmiiasympt}}
\begin{proof}
 To calculate the asymptotic for $t/\tau >> 1$, we apply Proposition \ref{Isampt} to Eq. \eqref{incvarfbmiieq}. For $t/\tau << 1$ asymptotic we transform the part of   the increment process of RL FBM containing the hypergeometric function into the Fox H function using formula (15.3.4) from \cite{Abramowitz} and 1.131 from \cite{Mathai2010}
        \begin{equation}\label{fromhypergeotoh}
            \begin{aligned}
                \frac{4H \left(t+\tau\right)^{H - \frac{1}{2}}t^{H + \frac{1}{2}}}{H + 1/2} \Hypergeometric{2}{1}{\frac{1}{2} - H, 1}{\frac{3}{2} + H}{\frac{t}{t + \tau}} = \frac{4H \tau^{H - \frac{1}{2}}t^{H + \frac{1}{2}}}{\Gamma\left(1/2-H\right)}\Hypergeometric{2}{1}{\frac{1}{2} - H, 1}{\frac{3}{2} + H}{-\frac{t}{\tau}} = \\ = \frac{4H \tau^{H - \frac{1}{2}}t^{H + \frac{1}{2}}}{\Gamma\left(1/2-H\right)}\FoxH{1}{2}{2}{2}{\frac{t}{\tau}}{H + \frac{1}{2}, \frac{1}{2} - H}{1, 1}{0, -H - \frac{1}{2}}{1, 1}.
            \end{aligned}
        \end{equation} Then, we get
    \begin{equation}
        \begin{aligned}
            \mathbb{E}\left[\left(b_{H}^{*\tau}\left(t\right)\right)^2 \right] = \left(t+\tau\right)^{2H} + t^{2H} - \frac{4H \tau^{H - \frac{1}{2}}t^{H + \frac{1}{2}}}{\Gamma\left(1/2-H\right)}\FoxH{1}{2}{2}{2}{\frac{t}{\tau}}{H + \frac{1}{2}, \frac{1}{2} - H}{1, 1}{0, -H - \frac{1}{2}}{1, 1} = \\ = \tau^{2H}\left(\left(1+\frac{t}{\tau}\right)^{2H} + \left(\frac{t}{\tau}\right)^{2H} - \frac{4H}{\Gamma\left(1/2-H\right)}\left(\frac{t}{\tau}\right)^{H+\frac{1}{2}}\FoxH{1}{2}{2}{2}{\frac{t}{\tau}}{H + \frac{1}{2}, \frac{1}{2} - H}{1, 1}{0, -H - \frac{1}{2}}{1, 1}\right) = \\ = \tau^{2H}\left(I_1 + I_2 - I_3\right).
        \end{aligned}
    \end{equation}
   To approximate $I_1$ we expand it in Taylor series.
    For $I_3$, when $H \neq \frac{1}{2}$, we apply formula 8.3.2.3 from \cite{Prudnikov} to represent H-function as the series. Thus, we have 
    \begin{equation}\label{i3series}
        \begin{aligned}
            I_3 = \frac{4H}{\Gamma\left(1/2-H\right)}\left(\frac{t}{\tau}\right)^{H+\frac{1}{2}}\FoxH{1}{2}{2}{2}{\frac{t}{\tau}}{H + \frac{1}{2}, \frac{1}{2} - H}{1, 1}{0, -H - \frac{1}{2}}{1, 1} 
            = \\ = \frac{4H}{\Gamma\left(1/2-H\right)}\left(\frac{t}{\tau}\right)^{H+\frac{1}{2}} \sum_{k=0}^{\infty} \frac{\Gamma\left(\frac{1}{2} - H + k\right)}{\frac{1}{2} + H + k}\frac{\left(-1\right)^k}{k!}\left(\frac{t}{\tau}\right)^{k} = .
        \end{aligned}
    \end{equation}
    Using \eqref{i3series} we obtain the asymptotic of $I_3$,
    \begin{equation}\label{i3asympt}
        \begin{aligned}
            I_3 \sim \frac{4H}{\Gamma\left(1/2-H\right)}\left(\frac{t}{\tau}\right)^{H+\frac{1}{2}} \left[\frac{\Gamma\left(\frac{1}{2} - H\right)}{\frac{1}{2} + H} - \frac{\Gamma\left(\frac{3}{2} - H\right)}{\frac{3}{2} + H}\frac{t}{\tau} + \frac{\Gamma\left(\frac{5}{2} - H\right)}{\frac{5}{2} + H}\frac{1}{2}\left(\frac{t}{\tau}\right)^{2}\right] = \\ = 4H\left(\frac{t}{\tau}\right)^{H+\frac{1}{2}} \left[\frac{1}{\frac{1}{2} + H} - \frac{\left(\frac{1}{2} - H\right)}{\frac{3}{2} + H}\frac{t}{\tau} + \frac{\left(\frac{1}{2} - H\right)\left(\frac{3}{2} - H\right)}{\frac{5}{2} + H}\frac{1}{2}\left(\frac{t}{\tau}\right)^{2}\right].
        \end{aligned}
    \end{equation}
    If $H=\frac{1}{2}$, then RL FBM reduces to the (ordinary) Brownian motion and $\mathbb{E}\left[\left(b_{H}^{*\tau}\left(t\right)\right)^2 \right] = \tau$. Combinig all the approximations done above and keeping only the dominating terms we obtain \eqref{structfuncfbmiiasympt}.
\end{proof}

\end{document}